\documentclass[12pts, oneside]{article}
\usepackage{geometry}
\usepackage{color}
\usepackage{amssymb}
\usepackage{color}
 \usepackage{graphicx}
\usepackage{multimedia}
\usepackage{amsmath,amssymb}
\usepackage{ulem}
\usepackage{soul}

\newtheorem{thm}{Theorem}[section]
 \newtheorem{cor}[thm]{Corollary}
 \newtheorem{lem}[thm]{Lemma}
 \newtheorem{prop}[thm]{Proposition}
 \newtheorem{defn}[thm]{Definition}
 \newtheorem{rem}[thm]{Remark}

\newcommand{\dvg}{{\rm div\,}}

\newcommand{\proof}{{\bf Proof\\}}

\begin{document}
\date{}

\title{Time-periodic solutions to  heated ferrofluid flow models}

\author{ Kamel Hamdache, \thanks
{ L\'eonard de Vinci P\^ole Universitaire, DVRC. 92 916 Paris La
D\'efense Cedex. France (Email: kamel.hamdache@devinci.fr)} \and
{Djamila Hamroun},  \thanks {Corresponding author}
\thanks{Laboratoire AMNEDP Facult\'e de Math\'ematiques,
Universit\'e USTHB, Algiers, Algeria (Email: djamroun@yahoo.fr)}
\and {Basma Jaffal-Mourtada}
\thanks { L\'eonard de Vinci P\^ole Universitaire, DVRC. 92 916
Paris La D\'efense Cedex. France (Email: basma.jaffal@devinci.fr)} }
\maketitle

{\bf Subjclass } {35Q35, 76D05}

{\bf Keywords } {Magnetic fluid flow, Navier-Stokes equations,  {convection-diffusion equation,  magnetostatic equations,} magnetization,  time-periodic weak  solutions.\\
}

\begin{abstract}
In this work we prove the existence of time-periodic solutions to a
model  describing  a ferrofluid flow  heated from below.
 Navier-Stokes equations satisfied by the fluid velocity are coupled to the temperature equation  and the magnetostatic
  equation satisfied by the magnetic potential. The magnetization is assumed to be parallel to the magnetic field and is given by a non-linear magnetization
   law generalizing the  Langevin law. The proof is based on  a semi-Galerkin approximation and regularization methods together with  the fixed point
    method.
\end{abstract}

\maketitle
\section{Introduction}

 The mathematical model describing heat transfer in an electrically
non-conducting incompressible ferrofluid heated from below is driven
by the fluid velocity ${ U}$, the pressure ${p}$, the temperature
$\tau$ and the magnetic field ${ H}$. The fluid velocity ${U}$ and
the pressure ${p}$ satisfy in a cylindrical domain  the
Navier-Stokes equations with the Kelvin force $\mu_0 ({
M}\cdot\nabla) {H}$  and we assume that the Oberbeck-Boussinesq
approximation is
 valid. Moreover the fluid is subjected to an external magnetic field and a gravity
 force.
  The temperature $\tau$ satisfies  a convection-diffusion equation
and the magnetic field $H$ satisfies
 the magnetostatic equations.
  \\
 The resulting magnetization ${ M}$ is assumed to be parallel to the
demagnetizing field ${H}$. Collinearity is a good approximation when
the internal rotation of colloidal particles can be neglected.
  Following the Langevin law \cite{ARRRP, PT}, the
magnetization writes as
\begin{equation}\label{langevin}
{M}= M_S{\mathcal L}\Big(\frac{\mu_f m\,{|{H}|}}{k_B\, \tau}
\Big)\frac{{ H}}{|H|},  \ \ \ 
{\mathcal L}(s)= \frac{1}{\tanh(s)}-\frac{1}{s},\,\,  s\geq 0,
\end{equation}
where    $|\cdot  |$ denotes the Euclidian norm in $\mathbb R^3$,
$M_S$ is the saturation magnetization, $m$ is the magnetic moment of
the particle and $k_B$ is the Boltzmann's constant.

In this work we shall use the following  general form of the
Langevin magnetization law
\begin{equation}\label{langevin2}
M= M_S\, (\tau_\star-\tau) \, \chi(|{ H}|) \frac{{ H}}{|{ H}|},
\end{equation}
where $0\leq \tau\leq \tau_\star$  and $\chi$ is assumed to be such
that $\chi(s)\geq 0$ for all $s> 0$
 and the function $\displaystyle \frac{\chi(s)}{s}$ admits a finite limit as $s\to 0$  in order to have the law (\ref{langevin2})
  well defined for all ${ H} $. Of course the Langevin function $\mathcal L$ verifies these two  conditions.
  In the sequel we will   make additional  hypotheses for $\chi$.
  The   example $\chi(s)= \arctan(b\, s)$
   with $b>0$  {which satisfies these assumptions} was used in \cite{OBM}.\\
  The Kelvin force represented by $\mu_0 (M\cdot\nabla) H$ where
$\mu_0$ is the magnetic permeability of vacuum, can be written as
\begin{equation} \nonumber
\mu_0 ({M}\cdot\nabla){ H}= \nabla p_m - \mu_0\, M_S\, \kappa(|H|)\,
\nabla (\tau_\star-\tau),
\end{equation}
where we set
 \begin{equation}\label{kappa}
\kappa(s) = \int_0^{s} \chi(r)\, dr, \ \ s\geq 0,
\end{equation}
and where the magnetic pressure is given by
\begin{equation}\nonumber
p_m  =\mu_0\, M_S \, (\tau_\star-\tau) \, \kappa(|H|).
\end{equation}

  Ferrofluid flows have been intensively studied by
many authors   for their interesting properties for engineering and
medical applications. In ~\cite{An, Ro, S2} for instance, the
authors investigate the stability analysis of solutions with respect
to various parameters. The Langevin magnetization law
(\ref{langevin}) has been considered recently in \cite{AH6} for a
stationary homogenization problem. A model of heated magnetic fluid
flow was studied in \cite{AH5}, from the point of view of the
existence of global weak solutions, using a linear approximation of
the Langevin magnetization law (\ref{langevin}).
 In this
paper,  a more  general problem  is considered,  with a generalized
(nonlinear) Langevin magnetization law,  impacting strongly the
force in the Navier-Stokes equations. For this new model, we aim
  to prove
existence of time-periodic solutions.
 Let us mention the work given
in \cite{OJ} where existence of time-periodic solutions are proved
for the Shliomis equation.

\subsection{The model equations}
We assume that   the ferrofluid occupies a cylindrical   domain
 $D= \Omega\times (0,
d)$ of height $d>0$ where the cross section $\Omega$ is a regular
open and bounded  domain
 of ${\mathbb R}^2$. The generic  point of $D$ is denoted $x =(\widehat{x},x_3)$  where    $\widehat{x}=(x_1,x_2)$ is the generic point of $\Omega$.
  The scalar product in   $\mathbb R^2$ and $\mathbb R^3$ 
    is denoted by $\cdot$ and the euclidian norm by $|\cdot
 |$. The operators $\nabla,\,  \dvg, \ \Delta$ are the usual differential
 operators in the $x$-variable and we denote by $\widehat\nabla, \
 \widehat\Delta$ the operators in the $\widehat x$-variable. \\
The lateral boundary of $D$ is denoted $\Sigma= \partial\Omega
\times (0,d)$ and the lower  and upper boundaries $\Gamma^-=
\Omega\times\{ 0 \}$ and $\Gamma^+=
\Omega\times\{  d \}$ respectively. We denote by $\nu$ the outward normal to the boundary of $D$  and  by $\widehat \nu$ the outward normal to the boundary $\partial \Omega$ of $\Omega$.\\
Let $T>0$, we set $\Omega_T=(0,T)\times \Omega$, $D_T=(0,T)\times
D$, $\Gamma_T^\pm= (0,T)\times \Gamma^\pm$ and  $\Sigma_T
=(0,T)\times \Sigma$.
\\
 The   lower boundary $\Gamma^- $ of $D$  is maintained at a given temperature {$\zeta=\zeta(t, \widehat{x})$}   which is $T$- periodic, that is
$\zeta(T)=\zeta(0)$ and we expect that the temperature $\tau$ is
less than
\begin{equation}\label{taustar}
\tau_\star=\|\zeta\|_{L^\infty(\Omega_T)}.
\end{equation}

\noindent The equations satisfied by  $({ U}, \tau, { H})$ and the
pressure $p$ in $D_T$ are, see~\cite{An, Ro, S2}
 \begin{equation}\label{systembis}
\begin{array}{ll}
\displaystyle {\rm div\,}{ U}=0, \\[2ex]
\displaystyle \rho_0 (\partial_t { U}+ ({U}\cdot\nabla){
U})-\mu\Delta
{ U}+\nabla {p} = {\mathcal S}, \\[2ex]
\displaystyle \rho_0 c_p(\partial_t {\tau}+ { U}\cdot\nabla{\tau}) -\eta \Delta \tau= 0 ,\\[2ex]
\displaystyle {H}=\nabla\varphi,\,\, \dvg\big(\nabla\varphi+  M_S\,
(\tau_\star-\tau) \, a(\nabla\varphi)\big)= F,
\end{array}
\end{equation}
where $\mu$ is the dynamical viscosity, $\eta$ is the thermal
conductivity coefficient of the fluid, $c_p$ is the specific heat at
constant pressure and $\rho_0$ is the reference density. \\
The force in the Navier-Stokes equation is given by
 \begin{equation}
\label{Sforce2bis} {\mathcal S}  = \mu_0\, M_S\, (\tau_\star-\tau)\,
\chi(|H|) \nabla |H|+\rho_0(1+ \alpha\, (\tau_\star-\tau))\, g,
\end{equation}
where the gravity force is such that  $g=-|g|\,  (0,0,1)$. In the
magnetostatic equations,  the source term ${F}$  is a given function
which is $T$- periodic, that is $F(T)=F(0)$ and the function
$a(\xi)$ in the flux is defined  by
\begin{equation}\label{flux}
a(\xi)=
 \chi(|\xi|)\, \frac{\xi}{|\xi|} \ \ \mbox{ if } \xi\in {\mathbb
R}^3\setminus \{0\}\ \ \mbox{ and }\ \  a(0)= 0.
\end{equation}
    System~\eqref{systembis}  is supplemented
with the boundary conditions
\begin{equation}\label{bcbis}
\begin{array}{ll}
\displaystyle {U}=0, \ \ \ \big(\nabla\varphi+  M_S\,
(\tau_\star-\tau) \, a(\nabla\varphi)\big)\cdot{\nu}=0\ \mbox{ on } \, (0,T)\times \partial D, \\[2ex]
 \displaystyle \tau = \zeta\, \mbox{ on } \, \Gamma_T^-,\ \ \ \ \eta\nabla
\tau\cdot {\nu}=0\,\, \mbox{ on } \, \Sigma_T \cup \Gamma_T^+,
\end{array}
\end{equation}
and with the time-periodic conditions
\begin{equation}\label{initperbis}
\begin{array}{ll}
\displaystyle { U}(0)={ U}(T),\,\, {\tau}(0)={\tau}(T)\ \ \ \mbox{
in } \, D.
\end{array}
\end{equation}
 The data $\zeta$   has to satisfy the compatibility condition $
\widehat\nabla \zeta\cdot\widehat \nu=0$ on $(0,T)\times
\partial \Omega$, which implies since $\zeta$ is independent of the vertical variable $x_3$, that  $  \nabla \zeta\cdot  \nu=0$ on
$\Sigma_T\cup \Gamma_T^+$. Other boundary conditions for the
temperature will be considered in the last section of this paper.
 Recall that  $(\tau_\star-\tau)\geq 0$ and $\chi(|\nabla \varphi|)\geq 0$, so the boundary condition satisfied by $\varphi$ reduces to
 the homogeneous Neumann one $\nabla\varphi\cdot \nu=0$ 
and since ${\varphi}$ is only determined up to a constant, we
additionally require that $\varphi$ is normalized, that is  $\int_D
\varphi(t,x)\, dx=0$ a.e. in $(0,T)$.
\\
It is useful to make the change of variable
  $ \widetilde \tau= \tau- \zeta$,
 rewriting the previous system of equations, we get
\begin{equation}\label{system}
\begin{array}{ll}
\displaystyle {\rm div\,}{ U}=0, \\[2ex]
\displaystyle \rho_0 (\partial_t { U}+ ({U}\cdot\nabla){
U})-\mu\Delta
{ U}+\nabla {p} = {\mathcal S}, \\[2ex]
\displaystyle \rho_0 c_p(\partial_t {\widetilde \tau}+ { U}\cdot\nabla{\widetilde \tau}) -\eta \Delta \widetilde \tau= -\rho_0 c_p \,  U \cdot\nabla \zeta+Z(\zeta),\\[2ex]
\displaystyle {H}=\nabla\varphi,\,\, \dvg\big(\nabla\varphi+  M_S\,
b(\widetilde \tau) \, a(\nabla\varphi)\big)= F,
\end{array}
\end{equation}
where we set
\begin{equation}\label{Z}
Z(\zeta)=-\rho_0 c_p \,\partial_t {\zeta}+ \eta \Delta \zeta,
\end{equation}
\begin{equation}\label{b}
b(\widetilde \tau)=   \tau_\star-\zeta-\widetilde \tau\geq 0.
\end{equation}
The force $\mathcal S$ and the boundary and time-periodic conditions
associated to (\ref{system}) read as
\begin{equation}
\label{Sforce2} {\mathcal S}={\mathcal S}(\widetilde \tau,|H|) =
\mu_0\, M_S\, b(\widetilde \tau)\, \chi(|H|) \nabla |H|+\rho_0(1+
\alpha\, b(\widetilde \tau))\, g,
\end{equation}
\begin{equation}\label{bc}
\begin{array}{c}
\displaystyle {U}=0, \ \ \ \big(\nabla\varphi+  M_S\,
b(\widetilde \tau) \, a(\nabla\varphi)\big)\cdot{\nu}=0\ \mbox{ on } \, (0,T)\times \partial D, \\[2ex]
\displaystyle \widetilde \tau = 0\, \mbox{ on } \, \Gamma_T^-,\ \ \
\ \eta\nabla \widetilde \tau\cdot  {\nu}=0\,\, \mbox{ on } \,
\Sigma_T\cup \Gamma_T^+,
\end{array}
\end{equation}
\begin{equation}\label{initper}
\begin{array}{ll}
\displaystyle { U}(0)={ U}(T),\,\, {\widetilde \tau}(0)={\widetilde
\tau}(T)  \ \ \mbox{ in } \, D.
\end{array}
\end{equation}

Problem \eqref{system}  with the boundary conditions~\eqref{bc} and
time-periodic conditions~\eqref{initper} is labeled problem
$({\mathcal P}_{per})$. Notice therefore that the magnetic potential
$\varphi$ and the magnetic field $H$ are necessarily $T$- periodic.
Of course a solution of problem $({\mathcal P}_{per})$ may be
extended by periodicity to the whole ${\mathbb R}$ and then
satisfies the property
\begin{equation}\label{propertyper}
U(t)= U(t+T),\,\, \widetilde\tau(t)=\widetilde\tau(t+T), \ \
H(t)=H(t+T), \ \ \,\,
\,\, t\in {\mathbb R}.
\end{equation}
In the sequel, we will also consider problem $({\mathcal P}(U_0,
\tau_0))$ (or $({\mathcal P})$ if no confusion arises) consisting of
the system of equations~\eqref{system}--\eqref{Sforce2} with the
boundary conditions~\eqref{bc} and   initial conditions
\begin{equation}\label{init}
U(0)= U_0 ,\,\, \widetilde\tau(0)= \tau_0.
\end{equation}

\subsection{Notations and main result}
$ {\bullet}$ {\bf Functional spaces.} Let $L^{q} $, $H^s $ and
$W^{s,q} $ ($1\le {q}\le \infty$, $s \in {\mathbb R}$) be the usual
Lebesgue and Sobolev spaces of scalar-valued functions, respectively
and we
  will use the notations $\mathbb L^{q} $, $\mathbb H^{s} $  and $\mathbb{W}^{s,q} $ when dealing with
  vector fields (and sometimes more general tensor fields).
  We denote
by $\|\cdot\|_q$ both the $L^q(D)$ and  $\mathbb L^q(D)$ norms and
simply by $\|\cdot\|$ if $q=2$. We will also make use of the well
known spaces $L^p(0,T; \mathcal V)$, $H^s(0,T; \mathcal V)$ and
$\mathcal
C([0,T]; \mathcal V)$ where $\mathcal V$ is a general Banach space.\\
 The duality product
between the Banach space
  ${\mathcal V}$ and its dual  ${\mathcal V}^\prime$ is denoted by $\langle\cdot ;
\cdot \rangle_{{\mathcal V}^\prime\times{\mathcal V}}$ or simply by
$\langle\cdot ; \cdot \rangle$ when there is no confusion of
notation and we denote by ${\mathcal C}([0,T]; {\mathcal V}^\prime
{-weak})$ the space of functions $v: [0, T] \longrightarrow {\mathcal V}^\prime$ which are continuous with respect to the weak topology.\\
Let $\mathcal D(D)$ the set of functions   which are infinitely
differentiable in $D$ with compact support in $D$ and  let $\mathbb
H^1_0 (D)$ the closure  of $(\mathcal
 D (D))^3$   in $\mathbb
H^1 (D)$. We introduce the classical functional spaces in the theory
of the Navier-Stokes equations, see~\cite{GA, GA2, LA, JLL, Tar,
TE},
$$\begin{array}{ll}
\displaystyle {\mathcal D} _{s} = \left\{ u\in (\mathcal
 D (D))^3;  \, {\rm div\,}u=0 \, \mbox{ in } D \right\},\\[2ex]\displaystyle
 {\mathcal U}_0 \, = \mbox{ closure of } {\mathcal D_s} \mbox{ in } \mathbb{L}^2(D), \\[2ex]\displaystyle
{\mathcal U}\, = \mbox{ closure of }{\mathcal D_s} \mbox{ in }
\mathbb{H}^1(D).\end{array}$$ We recall that
$$\begin{array}{c}
 {\mathcal U}_0= \left\{ u\in \mathbb{L}^2(D); \, {\rm div\,} u=0
\, \mbox{ in } D, \; u \cdot {\nu}=0  \, \mbox{ on }
\partial D \right\}, \\[2ex] {\mathcal U}=\left\{ u\in \mathbb{H}^1_0(D); \, {\rm div\,} u=0
\, \mbox{ in } D \right\}, \end{array}$$
and identifying ${\mathcal
U}_0$     with its dual, we get the inclusions ${\mathcal U} \subset
{\mathcal U}_0 \subset {\mathcal
U}\, ^\prime$.\\
 {For the temperature and magnetic potential, we
introduce the Hilbert spaces
\begin{align*}
 H^1_{\Gamma^-}&=\left\{ \theta\in  {H}^1(D); \,   \theta=0 \, \mbox{
on }     \Gamma^-\right\}, \ \ \ H^{-1}_{\Gamma^-}=(H^1_{\Gamma^-})^\prime, \\[1ex]  L^2_\sharp &= \{ w\in  L^2(D); \,   \int_D w\, dx=0
 \}, \ \ \ H^1_\sharp= L^2_\sharp \cap  {H}^1(D), \ \ \ H^{-1}_\sharp= (H^1_\sharp)^\prime.
\end{align*}
We recall the    Poincar\'e and Poincar\'e-Wirtinger inequalities:
there exists $C_p>0$ depending on $D$ such that
\begin{equation}\label{Poincare}
\|U\|\leq C_p \|\nabla U\|,  \ \ \forall U\in \mathcal U, \ \ \
 \|\theta\|\leq
C_p \|\nabla \theta\|,  \ \ \forall \theta
 \in H^1_{\Gamma^-},
 \end{equation}
 \begin{equation}\label{PW}
 \|w\|\leq C_p \|\nabla w\|, \ \ \forall w \in
 H^1_\sharp,
 \end{equation}
  so the spaces ${\mathcal U}$,
$H^1_{\Gamma^-}$  and  $ H^1_\sharp$  are   equipped with the $L^2-$
norm of the
gradient.\\
We will also need the following convex sets
\begin{equation}\label{tempspace}\begin{array}{c}
\mathcal W_0=\big\{\vartheta \in L^\infty(D), \ \ \  0 \leq
\vartheta+\zeta\leq \tau_\star\  \mbox{  a.e. in } D\big\},\\[2ex] \mathcal
W=\big\{\theta \in   \mathcal C([0,T];L^2(D) )\cap
L^\infty(D_T); 
0\leq \theta+\zeta\leq \tau_\star  \ \ \mbox{ a. e. in } D_T\big\}
.\end{array}
\end{equation}
\bigskip

\noindent $ {\bullet}$ {\bf Hypotheses.} All the physical parameters
are assumed to be nonnegative. In addition, we will make use of the
following hypotheses:
\begin{equation}\label{hyp}
\begin{array}{c}
F\in {\mathcal C}([0,T];L^2_\sharp)\cap L^\infty(0,T; L^3(D)),\ \ \
F(0)=F(T), \ \ \ 
g\in  \mathbb L^\infty(D_T), \ \ \  g(0)=g(T),
\end{array}
\end{equation}

 \begin{equation}\label{hyp
zeta}\begin{array}{c}\zeta\in {\mathcal C^1}([0,T];L^2(\Omega))\cap
\mathcal C([0,T]; H^2(\Omega) ), \ \ \ 
 \widehat\nabla
\zeta\cdot\widehat \nu=0 \ \mbox{ on } (0,T)\times
\partial \Omega, \\[2ex]  \zeta(0)=\zeta(T)\ \ \mbox{ in
}\Omega,  \zeta \geq 0 \ \ \mbox{a.e. in } \Omega_T, \ \
\tau_\star=\|\zeta\|_{L^\infty(\Omega_T)}>0,
\end{array}
\end{equation}
\smallskip

\begin{equation}\label{hyp chi}
\begin{array}{c}
\chi\in {\mathcal C}^1([0, +\infty[),  \ \ \chi(0)= 0, \ \ \  0\leq
\chi(s)\leq \chi_0,\ \ \   0\leq \chi^\prime(s)\leq \chi_1, \ \
\forall s\geq 0, \\[2ex] \mbox{ so } |\chi(r)-\chi(s)|\leq \chi_1|r-s|,
\ \ \forall \,  r, s\geq 0,
\end{array}
\end{equation}
where $\chi_0, \chi_1>0$ are constants. We deduce  that
\begin{equation}\label{hyp Z} Z(\zeta)=-\rho_0 c_p \,\partial_t {\zeta}+ \eta
\widehat\Delta \zeta \in {\mathcal C}([0,T];L^2(\Omega)),
\end{equation}
and the following properties for the function $\kappa$ given by
(\ref{kappa}) are satisfied
\begin{equation}\label{hyp kappa}
\begin{array}{c}
\kappa\in {\mathcal C}^2([0, +\infty[),\ \ \kappa(0)= 0, \ \ \kappa
\mbox{ is convex },\ \ \
  \kappa(s)\geq 0,  \\[2ex]|\kappa(r)-\kappa(s)|\leq
\chi_0|r-s|, \ \ \forall \,  r, s\geq 0. \end{array}
\end{equation}
\bigskip

\noindent  Before going on, let us check that we can find a function
$\zeta$ satisfying the properties given in (\ref{hyp zeta}).
 Indeed if we take a non-trivial function  $f\in \mathcal C^1([0,T]; L^2(\Omega))$  such that 
  $f(0)=f(T)$ and $\int_{\Omega_T}f\, dxdt=0$, then (see \cite{Bostan}, Theorem 3.12) there exits a function $\varsigma$ (which is unique up to a constant) such that \ $v\in   {\mathcal
C^1}([0,T];L^2(\Omega))\cap \mathcal C([0,T]; H^2(\Omega) )  $ \ and
satisfies equation
\begin{equation}\label{chalper}
\begin{array}{ll}
\rho_0c_p
\partial_t  \varsigma- \eta \widehat\Delta \varsigma =f \ \  \mbox{ in }\Omega_T, \ \  \\[2ex]
\widehat\nabla \varsigma\cdot\widehat \nu=0 \ \mbox{ on }
(0,T)\times
\partial \Omega, \ \ 
\varsigma(0)=\varsigma(T)\ \ \mbox{ in }\Omega .
\end{array}
\end{equation}
 Moreover since
  $H^2(\Omega)\subset
L^\infty(\Omega)$ so $\varsigma \in L^\infty(\Omega_T)$. Then
$\zeta=\varsigma+\|\varsigma\|_{L^\infty(\Omega_T)}$ verifies
(\ref{hyp zeta}).
 \\

With these notations at hand,  we can give the  definition of weak
solutions of our problem and formulate the main results of this
paper.
\begin{defn}\label{definition} Let $T>0$ be fixed and let  $F$, $g$, $\zeta$ and $\chi$ satisfying  (\ref{hyp}), (\ref{hyp zeta}) and (\ref{hyp
chi}). We say that $( U,  {\widetilde \tau}, {H})$ is a
time-periodic of period $T$ (or $T$-periodic) weak solution   of
problem $(\mathcal{P}_{per})$ if the conditions (i)--(vii) below are
satisfied:
\begin{enumerate}
\item[(i)] the velocity $U$ belongs to
 $L^\infty(0,T; {\mathcal U}_0)\, \cap\,  L^2(0,T; {\mathcal U}) \, \cap \, {\mathcal C}([0,T]; {\mathcal U}_0 \mbox{-weak})$\,\, {and}\,\,
$ U(0)=U(T)$,
\item[(ii)] the magnetic field $H$ is such that ${H}=\nabla \varphi$ where $\varphi$ belongs to $\mathcal C([0,T]; H^1_\sharp)$ $\cap\,  L^\infty(0,T;
W^{1,3}(D))$
and $\varphi(0)= \varphi(T)$,
\item[(iii)] the temperature $\widetilde \tau$ belongs to  ${\mathcal C}([0,T]; L^2(D)) \cap L^2(0,T; H^1_{\Gamma^-})$, $\widetilde \tau(0)= \widetilde \tau(T)$
 and  satisfies $0\leq \widetilde \tau+\zeta \leq  \tau_\star\; \mbox{ a.e. in }
 D_T$,
\item[(iv)]
the external body force ${\mathcal S}$ defined by (\ref{Sforce2})
belongs to $L^2(0,T;  {\mathcal U\, ^\prime})$ ,
\item[(v)] the momentum equation  holds
weakly, in the sense that, for every ${\Psi} \in {\mathcal U}$
\begin{equation}\label{eqU}
\begin{array}{r} \displaystyle
\rho_0\frac{d}{dt} \int_{D}{ U}\cdot{\Psi}\, dx +\rho_0 \int_{D} ({
U}\cdot\nabla){ U} \cdot{\Psi}\, dx + \mu \int_{D} \nabla{ U}\cdot
\nabla {\Psi}\, dx = \\[2ex] \displaystyle \langle {\mathcal S}, {\Psi} \rangle \ \mbox{
in } {\mathcal D}^\prime(]0, T[),\end{array}
\end{equation}
where  $\nabla { U}\cdot \nabla \Psi=\displaystyle\sum_{i=1}^3\nabla
{ U}_i\cdot \nabla \Psi_i$  and ${\mathcal D}^\prime(]0, T[)$ is the
space of distributions on $]0,T[$,
\item[(vi)] the magnetic potential $\varphi$ is a weak solution of the magnetostatic equation, i.e. for   $t\in [0,T]$,
\begin{equation}
\label{variat} \begin{array}{r}\displaystyle \int_D
\Big(\nabla\varphi(t)+ M_S\, b\big(\widetilde \tau(t)\big) \,
a\big(\nabla\varphi(t)\big)\Big) \cdot \nabla
\phi \, dx= \\[2ex] \displaystyle - \int_D F(t)\, \phi\, dx, \;\, \forall \phi \in
H^1_\sharp,\end{array}
\end{equation}
\item[(vii)] the temperature $\widetilde \tau$   satisfies
the integral identity
\begin{equation}
\begin{array}{lll}
\label{energyequality} && \displaystyle -\rho_0c_p\int_{D_T}
\widetilde \tau \big(\partial_t \theta + { U} \cdot \nabla
\theta\big) \, dxdt + \eta \int_{D_T} \nabla \widetilde \tau  \cdot
\nabla \theta \, dxdt
\\[2ex]
&& \displaystyle = -\rho_0c_p  \int_{D_T} ({ U}\cdot \nabla \zeta)
\, \theta \, dxdt+ \int_{D_T}Z(\zeta) \, \theta \, dxdt,
\end{array}
\end{equation}
for any ${\theta} \in C^1( \overline{D_T})$ such that
$\theta(0)=\theta(T)$ in $D$ and $\theta=0$ on ${\Gamma_T^-}$.
\end{enumerate}
\end{defn}

\bigskip

Our main result  is the following:
\begin{thm}\label{thm1}
Let $T>0$ be fixed. Assume hypotheses (\ref{hyp}),  (\ref{hyp zeta})
and (\ref{hyp chi}) hold true. Then there exists a \, $T$-periodic
solution $(U, \widetilde \tau, H=\nabla \varphi)$ of problem
$({\mathcal P}_{per})$ in the sense of  Definition
\ref{definition}.  Moreover there exists a pressure   $p\in
W^{-1,\infty}(0,T; L^2(D))$ such that the momentum equation is
satisfied in the sense of distributions.
\end{thm}

This result will be completed by  Theorem \ref{thm5} and  Theorem
\ref{thm4} given in the last section which are devoted  to existence
results of time-periodic weak solutions to similar models obtained
by considering other boundary conditions for the temperature $\tau$.
In Theorem \ref{thm5}, we consider the case of a
   magnetic fluid heated from below and above resulting in a Dirichlet
   condition on $\Gamma^-\cup \Gamma^+$,  then in Theorem \ref{thm4} we  consider the case  when the
   condition on $\Gamma^-$ is governed by the  Fourier law.

\medskip

  The plan of the paper is the following:\\
Section \ref{prelim} is devoted to  some preliminary results which
will be very helpful later and to
simplify the presentation,  the proofs will be done in an appendix.   \\
In Section \ref{Appro}, we present the regularization and
approximation methods on which the proof of our theorems is based.
We introduce the regularized problem $({\mathcal
P}_{per}^\varepsilon )$ of problem $({\mathcal P}_{per})$ depending
on a small parameter $\varepsilon>0$, obtained by regularizing the
intensity of the magnetic field in the force term $\mathcal S$ by a
smoothing sequence. A similar regularization is applied to the
initial and boundary problem $({\mathcal P})$ to get problem
$({\mathcal P}^\varepsilon )$.
 In solving the
regularized problems, we use a semi-Galerkin method   and we define
the approximated problems $({\mathcal P}_n^\varepsilon )$ and
$({\mathcal P}_{per, n}^\varepsilon )$ for $n\geq 1$. We rely on the
results of \cite{AH5} to get a solution of problem $({\mathcal
P}_n^\varepsilon )$ defined on $[0,T]$ and we
  prove its uniqueness. Then by using Schauder fixed
point theorem we prove the existence of a solution to problem
$({\mathcal P}_{per, n}^\varepsilon )$. \\ The next section is
devoted to obtain uniform estimates on the solutions
$(U_n^\varepsilon, \widetilde\tau_n^\varepsilon, H_n^\varepsilon)$
of problem $({\mathcal P}_{per, n}^\varepsilon )$, with respect to
$n$ allowing to perform the limit as $n\to + \infty$, $\varepsilon$
being fixed,  to get a solution $(U^\varepsilon,
\widetilde\tau^\varepsilon, H^\varepsilon)$ of problem $({\mathcal
P}_{per}^\varepsilon )$.\\  We achieve the proof of Theorem
\ref{thm1} in Section \ref{endproof}, by performing the limit as
$\varepsilon\to$ to get a solution of our
problem $({\mathcal P}_{per}  )$.\\
 Section \ref{other} is devoted to
the proofs of  the additional results Theorem \ref{thm5} and Theorem
\ref{thm4} on the models mentioned previously.
\bigskip

\begin{rem}\label{constants} In the sequel  we will denote by $C>0$ a
generic constant depending on the domain $D$,  the physical
constants of the problem, the data $T, F, g$, $\zeta$ and the bounds
$\chi_0, \chi_1$  related to $\chi$. If in addition the constant
depends on a special parameter like $\varepsilon>0$ or $n\in \mathbb
N$ then we will denote it  by $C_\varepsilon$ or $C_n$. If the
constant depends on various arguments $s_1, s_2, \cdots s_k $ then
we will denote it $C(s_1, s_2, \cdots s_k)$.
\end{rem}

 From now on, unless otherwise stated,
 we assume that   $F$, $g$, $\zeta$ and    $\chi$   satisfy hypotheses
(\ref{hyp}), (\ref{hyp zeta}) and (\ref{hyp chi}).

 \section{Preliminary results}\label{prelim}
 We give here some results   which will be very useful to us later. The proofs   are postponed to the appendix, see  Section \ref{appendix} .
  \subsection{The  magnetostatic equation}
Let  $\widetilde \tau$ be a given function  of the set $\mathcal W$
defined in (\ref{tempspace}). For $t\in [0,T]$, we consider the
magnetostatic equation
\begin{equation}\label{magperlim}
\begin{array}{ll}
\dvg\Big(\nabla\varphi+  M_S\, b(\widetilde \tau(t,.)) \,
a(\nabla\varphi)\Big)=
F(t,.)\ \  \mbox{in}\,\,   D,\\[2ex]\Big(\nabla\varphi+  M_S\,
b(\widetilde \tau(t,.)) \, a(\nabla\varphi)\Big)\cdot \nu=0 \ \
\mbox{on}\,\,
 \partial D,
\end{array}
\end{equation}
where $F$ satisfies hypothesis (\ref{hyp})   and $a( \xi)$ is
defined in (\ref{flux}). Notice that $b(\widetilde \tau)$ defined in
(\ref{b}) belongs to $ \mathcal C([0,T];L^2(D) )\cap \
L^\infty(D_T)$
 and satisfies $0\leq b(\widetilde \tau)\leq  \tau_\star$ in $D_T$.
\medskip

\noindent We have the following result:
\begin{prop}\label{prop1}
 Let $\widetilde\tau\in \mathcal W$, then
there exists a unique $\varphi\in \mathcal C([0,T]; H^1_\sharp )\,
\cap \, L^\infty(0,T;  W^{1,3}(D))$ satisfying  the magnetostatic
equation (\ref{magperlim}) for all $t\in [0,T]$ and
  $H=\nabla \varphi$ satisfies the estimates 
\begin{equation}\label{est mag}
\|H(t)\|\leq C\|F(t)\|,\ \ \ t\in [0,T],
\end{equation}
\begin{equation}\label{est magbis}
\|H(t)\|_3\leq C (1+\|F(t)\|_3), \ \ \ t\in (0,T).
\end{equation}
Moreover if $\widetilde\tau$ is $T$-periodic, then so is $\varphi $.
\end{prop}

\bigskip

Next we define the map
\begin{equation}\label{pot}
\mathcal H: \widetilde \tau \in  \mathcal W \longmapsto H= \nabla
\varphi,
\end{equation}
where $\varphi$ is the solution of the magnetostatic equation
(\ref{magperlim}) provided by Proposition \ref{prop1}.
 We prove  the following  result:

\begin{lem}\label{lem3}
It holds 
\begin{equation}\label{Hcont}
\|\mathcal H(\widetilde \tau_1)(t)-\mathcal H(\widetilde
\tau_2)(t)\|\leq C\, \|\widetilde \tau_1(t)-\widetilde \tau_2(t)\|,
\ \ \ \forall \widetilde \tau_1, \widetilde \tau_2\in \mathcal W , \
\ t\in [0,T].
\end{equation}
 In particular, if $(\widetilde \tau^m)_m\subset \mathcal W$ is such that  $\widetilde \tau^m \to  \widetilde \tau \ \ \mbox{strongly
 in }
   L^2(D_T)$ and $\widetilde \tau \in \mathcal W$ then $H^m=\mathcal H(\widetilde \tau^m) \to  H=\mathcal H(\widetilde \tau) \ \ \mbox{strongly in
   }
 \mathbb L^2(D_T)$.  Moreover  $H\in \mathcal C([0,T]; \mathbb L^2(D))\cap L^\infty (0,T; \mathbb L^3(D))$,    $H(t)=\nabla \varphi(t)$
  satisfies (\ref{magfaible}) for all $t\in [0,T]$ and  $H$ satisfies estimates (\ref{est
 mag}) and (\ref{est
 magbis}).
\end{lem}

\subsection{The temperature equation}\label{templemma}
 We consider the initial
boundary problem for the temperature equation
\begin{equation}\label{tempinit}
\begin{array}{ll}
\displaystyle \rho_0 c_p(\partial_t {\widetilde \tau}+ { U}\cdot\nabla{\widetilde \tau}) -\eta \Delta \widetilde \tau= -\rho_0 c_p \ { U}\cdot\nabla \zeta+Z(\zeta) \,\, \mbox{ in } \, D_T,\\[2ex]
\displaystyle \widetilde \tau =0\, \mbox{ on } \, \Gamma_T^-,\,\,
\eta\nabla
\widetilde \tau\cdot \nu =0\,\, \mbox{ on } \, \Sigma_T\cup \Gamma_T^+, \\[2ex]
\displaystyle {\widetilde \tau}(0)={\tau}_0 \,\, \mbox{ in } \, D,
\end{array}
\end{equation}
see (\ref{hyp Z}) for the definition of $Z(\zeta)$. We get the
following result:
\begin{prop}\label{temp}
Assume that $U\in
L^2(0,T; \mathcal U)$
 and
$\tau_0\in \mathcal W_0$.     Then  there exists a unique weak
solution $\widetilde \tau$ of  problem (\ref{tempinit}) such that
$$\widetilde \tau \in \mathcal C([0,T]; L^2(D))\cap L^2(0,
T;H^1_{\Gamma^-})\cap L^\infty(D_T),\ \ \  \partial_t \widetilde\tau
\in L^2(0,T; H^{-1}_{\Gamma^-})$$ which  satisfies
\begin{equation}\label{linfty} 0 \leq\widetilde \tau+\zeta\leq  \tau_\star
\mbox{ a.e. in} \ D_T,\end{equation} and
\begin{equation}\label{ener1bis}\|\widetilde \tau (t)\|^2+ \int_0^t \|\nabla
\widetilde \tau(s) \|^2\, ds \leq C(1+\| \tau_0 \|^2+ \int_0^t
\|U(s) \|^2\, ds ), \ \ t\in [0,T].\end{equation}
\end{prop}

  \bigskip

Hereafter, we introduce the operator
\begin{equation}\label{temperature}
\mathcal T: (U, \tau_0)\longmapsto \widetilde\tau,
\end{equation}
where $\widetilde\tau $ is the solution of problem (\ref{tempinit}).
 $\mathcal T$ satisfies the following properties:
\begin{lem}\label{taucont} Let $(U^i, \tau_0^i)\in  
L^2(0,T; \mathcal U)\times {\mathcal W}_0$ and let
$\widetilde\tau^i= \mathcal T(U^i, \tau_0^i)$ for $i=1,2$. Then it
holds for $t\in [0,T]$
\begin{equation}\label{contau}\begin{array}{r}\displaystyle \|\widetilde\tau^1(t)-\widetilde\tau^2 (t)\|^2+ \int_0^t \|\nabla
\widetilde\tau^1(s)-\nabla \widetilde\tau^2(s) \|^2\, ds \leq 
 C(\|\tau_0^1-\tau_0^2
\|^2+ \int_0^t \|U^1(s)-U^2(s) \|^2\, ds ).\end{array}\end{equation}
In particular if $(U^m, \tau_0^m)_m\subset 
L^2(0,T; \mathcal U)\times {\mathcal W}_0$ is such that $(U^m,
\tau_0^m) \to (U, \tau_0)$ strongly in  $L^2(0,T; \mathcal U_0
\times L^2(D))$ and $(U,\tau_0)\in 
L^2(0,T; \mathcal U)\times {\mathcal W}_0  $   then we have
$$\widetilde\tau^m=\mathcal T(U^m, \tau_0^m) \to \widetilde\tau=\mathcal T(U,\tau_0) \ \
\mbox{
  strongly  in }  L^2(D_T),$$ and $\widetilde\tau$ satisfies  problem (\ref{tempinit}) and estimates (\ref{linfty}) and
(\ref{ener1bis}).
 \end{lem}

 \subsection{ The regularized force $\mathcal
S_\varepsilon$}\label{Approforce}

In Section \ref{Appro}, we will consider the Navier-Stokes equation
with the force $\mathcal S$ replaced by the following one
\begin{equation}
\label{Sforcebis2} {\mathcal S}_\varepsilon ={\mathcal
S}_\varepsilon (\widetilde \tau,|H|) = \mu_0\, M_S\, b(\widetilde
\tau) \, \chi(\sigma^\varepsilon\star|H|)
\nabla(\sigma^\varepsilon\star|H|)+\rho_0(1+ \alpha\, b(\widetilde
\tau))\, g,
\end{equation}
where $\varepsilon>0$ is a small parameter and
$(\sigma^\varepsilon(x)) \subset \mathcal D(\mathbb R^3)$ is a
smoothing function   verifying
  $\sigma^\varepsilon\geq 0$ and  $\int_{{\mathbb R}^3} \sigma^\varepsilon(x)\, dx =1$.
   For a function $h$ defined on $D$, we set $\sigma^\varepsilon\star  h(x)=
  \int_{{\mathbb R}^3}\sigma^\varepsilon(y)  \overline h(x-y)\, dy $ where   $\overline h$ is the  extension of $h$ by 0 outside
  $D$.\\

 \noindent  
 We will establish the following properties:
    \begin{lem}\label{Sprop}
  Let $\widetilde \tau\in \mathcal W$ and $H= \mathcal H (\widetilde \tau)$
  then ${\mathcal
S}_\varepsilon= {\mathcal S}_\varepsilon (\widetilde \tau,|H|) \in
\mathcal C([0,T]; \mathbb L^2(D))$ with
\begin{equation}\label{Sbound1}
\|{\mathcal S}_\varepsilon (t) \|\leq C_\varepsilon,  \ \ t\in
[0,T],
\end{equation}
and if for $k=1,2$, $\widetilde \tau^k\in \mathcal W$ and $H^k=
\mathcal H (\widetilde \tau^k)$
 then  the forces ${\mathcal S}^k_\varepsilon= {\mathcal S}
_\varepsilon(\widetilde \tau^k, |H^k|)$ verify
\begin{equation}\label{Scont}
\|{\mathcal S}^1_\varepsilon(t)-{\mathcal S}^2_\varepsilon(t)\|
 \leq C_\varepsilon \|\widetilde \tau^1(t)-\widetilde \tau^2(t)\|, \ \ \ t\in [0,T].
\end{equation}
In addition if  $\widetilde \tau \in \mathcal W\cap L^2(0,T;
H^1_{\Gamma^-})$ then ${\mathcal S}_\varepsilon \in L^2(0,T;
{\mathcal U^\prime})$ and it holds
\begin{equation}\label{Sbound2}
\|{\mathcal S}_\varepsilon (t)\|_{\mathcal U^\prime}\leq C\, \big(1
+  \|\nabla\widetilde \tau (t) \|  \big), \ \ \ t\in (0,T).
\end{equation}
  \end{lem}

\section{ The regularized and approximated problems }\label{Appro}
  Let $\varepsilon>0$ be a small parameter, we will use the notations of Subsection \ref{Approforce}.\\   We
  define the regularized
  problems
$({\mathcal P}^\varepsilon)$ and  $({\mathcal P}^\varepsilon_{per}
)$ as  problems $({\mathcal P})$ and  $({\mathcal P_{per}})$ where
we replace  the force ${\mathcal S}$   by the regularized force
${\mathcal S}_\varepsilon$ defined in (\ref{Sforcebis2}).\\
We shall prove the following result:
\begin{prop}\label{thm3}
Under  hypotheses (\ref{hyp}), (\ref{hyp zeta}) and (\ref{hyp chi}),
there exists a $T$-periodic weak solution $(U^\varepsilon,
\widetilde\tau^\varepsilon, H^\varepsilon=\nabla
\varphi^\varepsilon)$ of problem $({\mathcal P}_{per}^\varepsilon)$
such that
$$\begin{array}{c}
U^\varepsilon\in L^\infty(0,T; \mathcal U_0)\cap  L^2(0,T; \mathcal
U)\cap {\mathcal
C}([0,T]; \mathcal U_0 \mbox{-weak}),\\[2ex]
\widetilde \tau^\varepsilon\in \mathcal C([0,T]; L^2(D))\cap
L^2(0,T;
H^1_{\Gamma^-}), \\[2ex]  H^\varepsilon=\nabla \varphi\in {\mathcal
C}([0,T]; {\mathbb L}^2(D))\cap L^\infty(0,T; {\mathbb
L}^3(D)).\end{array}
$$
Moreover it holds
\begin{equation}\label{bound9}
\begin{array}{cl}
0\leq \widetilde \tau^\varepsilon+\zeta\leq  \tau_\star \,\, \mbox{ a.e in } D_T,\\[3ex]
\|U^\varepsilon\|_{L^\infty(0,T; \mathcal U_0 )\cap L^2(0,T;
{\mathcal U})} +\|H^\varepsilon\|_{L^\infty(0,T; \mathbb
L^3(D))}\leq C,
\\[3ex]
\|\widetilde \tau^\varepsilon\|_{L^\infty(0,T; L^2(D))\cap L^2(0,T; H^1_{\Gamma^-})}+ \|\partial_t \widetilde \tau^\varepsilon\|_{L^2(0,T;  H^{-1}_{\Gamma^-} )}\leq C,\\[3ex]
 \|\mathcal S_\varepsilon(\widetilde \tau^\varepsilon, |H^\varepsilon|)\|_{L^2(0,T; {\mathcal
 U}^\prime)}+\|U^\varepsilon\otimes U^\varepsilon\|_{L^\frac{4}{3}(0,T;\mathbb L^{2}(D))}\leq C,
\end{array}
\end{equation}
\begin{equation}\label{bound9bis}
\|\partial_t U^\varepsilon\|_{L^\frac{4}{3}(0,T; {\mathcal
U}^\prime)}\leq C,
\end{equation}
where for $V,W\in \mathbb R^3$, \  $V\otimes W= (V_iW_j)_{ij}$.
\end{prop}

\medskip

\noindent {\bf Sketch of the proof. } 
 To begin with, we consider problem $({\mathcal P}^\varepsilon
(U_0, \tau_0) )$ (or shortly $({\mathcal P}^\varepsilon )$)
associated with the initial data
\begin{equation}\label{Hyp CI}
(U_0, \tau_0) \in \mathcal U_0\times \mathcal W_0,
\end{equation}
and  as in \cite{AH5},
 we  use a  semi-Galerkin method  to solve it.
 We consider  an orthogonal  basis  $(A_j)_{j\geq 1}\subset{\mathcal
U}\cap {\mathbb H}^2(D)$  of ${\mathcal U}$ which is  orthonormal in
${\mathcal U_0}$ and  for $n\geq 1$, we denote by $\mathcal V^n$ the
subspace spanned by $A_1, A_2, \cdots, A_n$. A function $ V_n=
\sum_{j=1}^n a_j^n A_j \in \mathcal V^n$ may be identified with the
vector
$a^n \in \mathbb R^n$ and we observe that    $ \|V_n\|=  |a^n|$ where $|\cdot |$ is the Euclidian norm in $\mathbb R^n$.  \\

We look for approximated solutions $ (U^\varepsilon_n,
\widetilde\tau^\varepsilon_n, H^\varepsilon_n=\nabla
\varphi_n^\varepsilon)$ of problem $({\mathcal P}^\varepsilon )$,
denoted simply $(U_n,\widetilde\tau_n, H_n)$ by omitting the
superscript $\varepsilon$,
proceeding   as follows:\\
 Let for $n\geq 1$, 
$  U_{0n}= \sum_{j=1}^n \alpha_{0j}^n A_j\in \mathcal
V^n$ , ($\alpha_{0}^n\in \mathbb R^n) $ 
be  the approximation of the  initial data $U_0$ and we assume that
as $n\to \infty$,
 $U_{0n} \to U_0 \ \ \mbox{
strongly in } \mathcal U_0. $ Then we look for $U_n$ of the form
$U_n(t)= \sum_{j=1}^n \alpha_j^n (t)A_j$,
  satisfying  the nonlinear system of differential equations
\begin{equation}\label{systemn1}
\begin{array}{rr}
\displaystyle \rho_0 \frac{d}{dt}\int_D U_n\cdot A_j\, dx +\rho_0
\int_D ( U_n\cdot\nabla )U_n\cdot A_j \, dx + \\[3ex] \displaystyle \mu \int_D\nabla
U_n\cdot \nabla A_j\, dx
=  \int_D {\mathcal S}_\varepsilon^{n}  \cdot A_j\, dx, \ \ \ j=1,\cdots, n,\\[3ex]
U_n(0)= U_{0n},
\end{array}
\end{equation}
where (see (\ref{Sforcebis2}), (\ref{temperature}) and (\ref{pot})
for the definitions of the operators ${\mathcal S}_\varepsilon$,
$\mathcal T$ and $\mathcal H$)
\begin{equation}\label{Sn}\begin{array}{c}
{\mathcal S}_\varepsilon^{n}= {\mathcal S}_\varepsilon (\widetilde
\tau_n, |H_n|)
  = \mu_0\, M_S  \,
 b(\widetilde \tau_n )\,\chi(\sigma^\varepsilon\star|H_n|)
\nabla(\sigma^\varepsilon\star|H_n|)  + 
\rho_0(1+ \alpha\, b(\widetilde \tau_n
))\, g,\end{array}
\end{equation}
  $$\widetilde \tau_n=\mathcal T(U_n,   \tau_0), \ \ \  H_n=\nabla \varphi_n=
\mathcal H(\widetilde \tau_n) .$$  So  the following equations are
satisfied
\begin{equation}\label{systemn2}
\begin{array}{ll}
\displaystyle \rho_0c_p(\partial_t \widetilde \tau_n +  U_n\cdot\nabla\widetilde \tau_n) -\eta \Delta\widetilde \tau_n  =-\rho_0c_p \   U_n\cdot \nabla \zeta+Z(\zeta),\\[2ex]
\widetilde \tau_n =0\ \ \mbox{on}\  \Gamma_T^-,\ \   \eta\nabla
\widetilde \tau_n\cdot \nu
=0\ \ \mbox{on}\  \Sigma_T\cup \Gamma_T^+,\\[2ex]\widetilde \tau_n(0)=  \tau_{0}\ \ \mbox{in}\  D,
\end{array}
\end{equation}
and for $t\in [0,T]$
\begin{equation}\label{magn}
\begin{array}{ll}\dvg(\nabla \varphi_n+  M_S\  b(\widetilde \tau_n) \
a(\nabla\varphi_n))= F(t)\ \  \mbox{in}\   D, \\[2ex] (\nabla \varphi_n+
M_S\  b(\widetilde \tau_n) \  a(\nabla\varphi_n))\cdot \nu=0\ \
\mbox{on}\ \
\partial D.
\end{array}
\end{equation}
 System (\ref{systemn1})-(\ref{systemn2})-(\ref{magn}) is labeled  problem $(\mathcal
 P_{n}^\varepsilon (U_{0n},  \tau_0))$ (or shortly $(\mathcal
 P_{n}^\varepsilon)$ if no confusion arises). Replacing in equations (\ref{systemn1}) and (\ref{systemn2}) the
initial conditions by the periodic ones, we obtain problem
$(\mathcal P_{per,n}^\varepsilon)$.\\

 Existence of a solution to problem $(\mathcal
P_{per,n}^\varepsilon)$ will be proved according to the following
steps.\\ First we consider the solution $(U_n,\widetilde\tau_n,
H_n=\nabla \varphi_n)$ of problem $(\mathcal
 P_{n}^\varepsilon (U_{0n}, \tau_0))$ obtained following the results given in \cite{AH5} by using a fixed point method and,  to be complete, we
  give the basic estimates which lead to this existence result.\\
 Then,   we
 prove that such a solution is unique allowing to introduce the operator
 ${\mathcal K}: (U_{0n}, \tau_0)\mapsto {\mathcal K}(U_{0n}, \tau_0)= (U_n(T),\widetilde\tau_n(T)) $. In proving
 that ${\mathcal K}$ possesses a fixed point by using Schauder's fixed point theorem, we get a solution of problem
$(\mathcal P_{per,n}^\varepsilon)$.\\
 Finally letting $n\to \infty$,
we get a solution of problem $(\mathcal
P_{per}^\varepsilon)$.\\

\subsection{Solving problem $(\mathcal P_ n^\varepsilon (U_{0n}, \tau_0))$}
Let $\varepsilon >0$ and $n\geq 1$ be fixed. We have  the following
result:
\begin{prop}\label{prop2}
Let $(U_{0n}, \tau_0)\in \mathcal V^n\times \mathcal W_0$. Under
hypotheses (\ref{hyp}), (\ref{hyp zeta}) and (\ref{hyp chi}), there
exists a solution $(U_n , \widetilde \tau_n , H_n =\nabla \varphi_n
)$ of problem $({\mathcal P}_{n}^\varepsilon  (U_{0n},  \tau_0))$
defined on $[0,T]$ such that
$$\begin{array}{c}U_n \in {\mathcal C}([0,T]; {\mathcal U}),\ \
\widetilde \tau_n \in \mathcal C([0,T]; L^2(D))\cap L^2(0,T;
H^1_{\Gamma^-}), \\[2ex]  H_n =\nabla \varphi_n\in {\mathcal
C}([0,T]; {\mathbb L}^2(D))\cap L^\infty(0,T; {\mathbb
L}^3(D)),\end{array}
$$
satisfying
\begin{equation}\label{b1} 0\leq \widetilde \tau_n+\zeta\leq  \tau_\star \ \ \mbox{a.e. in }
  D_T,\end{equation} and  for $ t\in [0,T]$
\begin{equation}\label{bound2}
\begin{array}{ll}
\displaystyle \|U_n(t)\|^2  +  \int_0^t \|\nabla U_n(s)\|^2 \, ds
\leq C\, \|U_{0n}\|^2+ C_\varepsilon,
\end{array}
\end{equation}
\begin{equation}\label{bound11} \|\widetilde \tau_n (t)\|^2+
\int_0^t \|\nabla \widetilde \tau_n(s) \|^2\, ds \leq C\, (1+\|
\tau_0 \|^2+ \|U_{0n}\|^2)+ C_\varepsilon,
\end{equation}
\begin{equation}\label{est mag n}
\displaystyle \|H_n(t)\|_3  \leq   C,
\end{equation}
\begin{equation}\label{boundlinfty}
\displaystyle \|U_n(t)\|_\infty  \leq  C_n (\|U_{0n}\|+
C_\varepsilon ).
\end{equation}
 Moreover this solution is unique.
\end{prop}
 \proof   To prove the existence of $(U_n,\widetilde
\tau_n, H_n)$, we can follow the method used in \cite{AH5} with  a
few precautions, given that in  our problem, the Langevin
magnetization law is nonlinear and the Dirichlet boundary condition
for the temperature is not constant.  Using the preliminary results
of Section \ref{prelim},   we  obtain a local solution defined on a
time interval $[0,T_n]$, $0<T_n\leq T$, such that
$$\begin{array}{c}U_n\in {\mathcal C}([0,T_n]; {\mathcal U}),\ \
\widetilde \tau_n\in \mathcal C([0,T_n]; L^2(D))\cap L^2(0,T_n;
H^1_{\Gamma^-}), \\[2ex]  H_n=\nabla \varphi_n\in {\mathcal C}([0,T_n];
{\mathbb L}^2(D))\, \cap \ L^\infty(0,T_n; {\mathbb
L}^3(D)).\end{array} $$ Let us prove that this solution can be
continued on $[0,T]$.  First in view of the results given in
Proposition \ref{temp}, the temperature $\widetilde \tau_n$
satisfies
$$0\leq \widetilde \tau_n+\zeta\leq  \tau_\star \ \ \mbox{a.e. in }
 (0,T_n)\times D,$$
\begin{equation}\label{bound1}
\begin{array}{r}
 \displaystyle \|\widetilde \tau_n (t)\|^2+ \int_0^t \|\nabla \widetilde \tau_n(s) \|^2\, ds
 \leq 
C(1+\| \tau_0 \|^2+ \int_0^t \|U_n(s) \|^2\, ds ), \ \ t\in [0,T_n].
\end{array}
\end{equation}
  Regarding  velocity $U_n$, we have
$$
\displaystyle\frac{\rho_0}{2}\frac{d}{dt} \|U_n\|^2+  {\mu} \|\nabla
U_n\|^2=\int_D {\mathcal S}^n_\varepsilon \cdot U_n\, dx,
$$
where ${\mathcal S}^n_\varepsilon$ is given by (\ref{Sn}), so using
Lemma \ref{Sprop} and Poincar\'e inequality, we get
 $$
\displaystyle\frac{\rho_0}{2}\frac{d}{dt} \|U_n\|^2+  {\mu} \|\nabla
U_n\|^2\leq C_\varepsilon \|U_n\|\leq C_\varepsilon  + \frac{\mu}{2}
\|\nabla U_n\|^2,
$$
 that is
 \begin{equation}\label{b2}
\begin{array}{ll}
\displaystyle\frac{\rho_0}{2}\frac{d}{dt} \|U_n\|^2 + \frac{\mu}{2}
\|\nabla U_n\|^2 \leq C_\varepsilon. \end{array}
\end{equation}
Finally we obtain   for $ t\in [0,T_n]$
$$
\begin{array}{ll}
\displaystyle \|U_n(t)\|^2  +  \int_0^t \|\nabla U_n(s)\|^2 \, ds
\leq C\, \|U_{0n}\|^2+ C_\varepsilon t\leq C\, \|U_{0n}\|^2+
C_\varepsilon T,
\end{array}
$$
and coming back to (\ref{bound1}), we deduce that for $ t\in
[0,T_n]$ $$
\|\widetilde \tau_n (t)\|^2+ \int_0^t \|\nabla \widetilde \tau_n(s)
\|^2\, ds \leq C(1+\| \tau_0 \|^2+ T\|U_{0n}\|^2)+ C_\varepsilon
T^2.
$$
The constants $C $  and $C_\varepsilon$
 are independent of $n$ so we
conclude that we can extend the solution $(U_n,\widetilde
\tau_n,H_n)$ to the whole interval $[0,T]$ and these bounds are
satisfied on $[0,T]$. (\ref{est mag n}) is a direct consequence of
(\ref{est mag}) and to prove (\ref{boundlinfty}), recall that $U_n$
is of the form $U_n= \sum_{j=1}^{n} \alpha_j^{n} A_j $ so
$$\|
U_n(t) \|_\infty \leq  \sum_{j=1}^n |\alpha_{j}^{n}(t) |\
\|A_j\|_\infty \leq C_n |\alpha^{n}(t)| \leq C_n \| U^n(t) \|,$$
   then by inequality
(\ref{bound2}), we deduce that
$$\|
U^n (t)\|_\infty \leq C_n (\|U_{0n}\|+ C_\varepsilon ).$$ The
uniqueness will follow from the lemma below. \hfill $\Box$

\medskip

\begin{lem}\label{lem4}
Let for $k=1,2$, $(U_{0n}^k,  \tau_{0}^k)\in  \mathcal V^n \times
\mathcal W_0$ and  $(U^k_n, \ \widetilde \tau^k_n, H^k_n=\nabla
\varphi_n^k)$ a solution of problem $({\mathcal P}_n^\varepsilon
(U_{0n}^k, \tau_{0}^k))$ provided by Proposition \ref{prop2}.
  Then there exists a constant $ C({n,\varepsilon}, \|U_{0n}^2\|)>0$ such that we
  have
\begin{equation}\label{uniq}
\begin{array}{r}
 \|U_n^1(t)-U_n^2(t)\|^2+  \|\widetilde \tau_n^1(t)-\widetilde \tau_n^2(t)\|^2 +\|H_n^1(t)-H_n^2(t)\|  \leq  \\[3ex]  C e^{t\, C({n,\varepsilon}, \|U_{0n}^2\|)}\, \big(\|U_{0n}^1-U_{0n}^2\|^2+
   \|  \tau_{0}^1-  \tau_{0}^2\|^2 \big), \ \ \ t\in [0,T],
\end{array}
\end{equation}
\begin{equation}\label{uniqbis}
\begin{array}{r}
\displaystyle  \int_0^T(\|\nabla U_n^1(t)-\nabla U_n^2(t)\|^2+
\|\nabla\widetilde \tau_n^1(t)-\nabla\widetilde \tau_n^2(t)\|^2)\,
dt   \leq
 \\[3ex]C e^{T\, C({n,\varepsilon}, \|U_{0n}^2\|)}\, \big(\|U_{0n}^1-U_{0n}^2\|^2+
   \|  \tau_{0}^1-  \tau_{0}^2\|^2 \big).
\end{array}
\end{equation}
 In particular
taking $U_{0n}^1= U_{0n}^2=U_{0n}$ and $  \tau_{0}^1= \tau_{0}^2=
\tau_0$ we deduce that  the solution $(U_n, \widetilde \tau_n, H_n)$
of problem $({\mathcal P}_n^\varepsilon (U_{0n},  \tau_{0}))$
  is unique.
\end{lem}
\proof We set  $U_n= U_n^1-U_n^2$, $\widetilde \tau_n= \widetilde \tau_n^1-\widetilde \tau_n^2$, 
$U_{0n}= U_{0n}^1-U_{0n}^2$, $  \tau_{0n}=  \tau_{0}^1-  \tau_{0}^2$
and for simplicity, we omit the subscript $n$. Then $U$ satisfies
the following system of equations
 $$\begin{array}{cl}
\displaystyle \rho_0 \frac{d}{dt}\int_D U\cdot A_j\, dx  +\mu
\int_D\nabla U\cdot \nabla A_j\, dx = \mathcal
F^1_j- \mathcal F^2_j, \ \ j=1, \cdots , n\\[2ex]
\displaystyle U(0)=U_{0},
\end{array}
$$
where for $j=1, \cdots , n$  and $k=1,2$
 $$
\mathcal F^k_{j}=-\rho_0 \int_D(  U^k\cdot\nabla )U^k \cdot A_j \,
dx +\int_D {\mathcal S}^k _\varepsilon\cdot A_j\, dx, \ \ \
{\mathcal S}^k _\varepsilon={\mathcal S} _\varepsilon(\widetilde
\tau^k, |H^k|).
$$
 Equation of $\widetilde \tau$ reads as
$$
\begin{array}{l}
 \rho_0c_p(\partial_t \widetilde \tau +  U^1\cdot\nabla\widetilde \tau ) -\eta \Delta
\widetilde \tau  = -\rho_0c_p\,
 U\cdot \nabla (\widetilde \tau ^2+\zeta )\ \mbox{ in }\  D_T,\\[2ex]
   \widetilde \tau =0\, \mbox{ on }\,
\Gamma^-_T,\ \  \eta\nabla\widetilde \tau \cdot \nu =0\ \mbox{ on }\
\Sigma_T\cup \Gamma^+_T ,\\[2ex]
\widetilde \tau (0)=  \tau_0\ \mbox{ in }\  D,
\end{array}
$$
so $\widetilde \tau$  satisfies the following inequality, see
(\ref{contau1}) in proof of Lemma \ref{taucont}
\begin{equation}\label{cont2}
 {\rho_0c_p} \frac{d}{dt}\|\widetilde \tau\|^2+ \eta \|\nabla \widetilde \tau\|^2\leq  C\|
 U\|^2,
\end{equation}
and    $U$ satisfies
$$\begin{array}{r}\displaystyle
\frac{\rho_0}{2}\frac{d}{dt} \|U\|^2 +\mu  \|\nabla U\|^2 =-\rho_0
\int_D \big(  (U^1\cdot\nabla )U^1 -(U^2\cdot\nabla )U^2\big) \cdot
U\, dx+ 
\int_D ( {\mathcal S}^1_\varepsilon-{\mathcal
S}^2_\varepsilon)\cdot U\, dx,
\end{array} $$
where by using (\ref{Scont})
\begin{equation}\label{Scont2}
\int_D (\mathcal S^1_\varepsilon- \mathcal S^2_\varepsilon) \cdot  U
\, dx\leq C_\varepsilon \|\widetilde \tau\| \,\| U\| \leq
C_\varepsilon \, (\|\widetilde \tau\|^2+\| U\|^2) .
\end{equation}
Now since
$$(U^1\cdot\nabla )U^1 -(U^2\cdot\nabla )U^2= (U^1\cdot\nabla )U +(U\cdot\nabla )U^2,$$
and $$\int_D (U^1\cdot\nabla )U \cdot U\, dx=0,  \ \  \  \ \ \int_D
(U\cdot\nabla )U^2 \cdot U\, dx= -\int_D (U\cdot\nabla )U \cdot
U^2\, dx,
$$
then
$$\begin{array}{c}\displaystyle |\int_D \big(  (U^1\cdot\nabla )U^1 -(U^2\cdot\nabla )U^2\big) \cdot
U\, dx|\leq \|U\|\ \|\nabla U \|\ \| U^2 \|_\infty.
 \end{array}$$ Using  inequality (\ref{boundlinfty})
for $U^2$, we deduce that
$$\begin{array}{r}\displaystyle|\int_D \big(  (U^1\cdot\nabla )U^1 -(U^2\cdot\nabla )U^2\big) \cdot
U\, dx| \leq C_n (\| U^2_{0n} \|+C_\varepsilon )\, \|U\|\ \|\nabla U
\| \\[2ex]\displaystyle \leq  \frac{\mu}{2}\|\nabla U \|^2 + \frac{1}{2\mu}
C_n^2 (\| U^2_{0n} \|+C_\varepsilon )^2 \, \|U\|^2,\end{array} $$
  which together with    inequality (\ref{Scont2}) lead   to
\begin{equation}\label{cont1}\displaystyle
 {\rho_0} \frac{d}{dt} \|U\|^2 +\mu  \|\nabla U\|^2\leq \big(\frac{1}{\mu} C_n^2 (\| U^2_{0n}
\|+C_\varepsilon )^2+2C_{ \varepsilon}\big) \ (\|U\|^2+ \|\widetilde
\tau\|^2).
\end{equation}
\noindent Finally summing up inequalities (\ref{cont2}) and
(\ref{cont1}), we can write
\begin{equation}\label{cont3}
\begin{array}{r}\displaystyle
   \rho_0 \frac{d}{dt}(\| U\|^2+c_p
 \|\widetilde \tau\|^2 )+\mu \|\nabla U\|^2+\eta \|\nabla \widetilde \tau\|^2  \leq  
 C({n,\varepsilon}, \|U_{0n}^2\|) \  \rho_0(\|U\|^2+c_p \|\widetilde \tau\|^2),
\end{array}
\end{equation} where $C({n,\varepsilon}, \|U_{0n}^2\|)>0$    depends  in an increasing way on $\| U^2_{0n}
\|$. Therefore by Gronwall  lemma, we deduce for all $t\in [0,T]$
$$\| U(t)\|^2+c_p
 \|\widetilde \tau(t)\|^2 \leq e^{t\, C({n,\varepsilon}, \|U_{0n}^2\|)}\ (\| U_0\|^2+c_p
 \| \tau_0\|^2),  $$
and using  (\ref{Hcont}), we arrive to (\ref{uniq}).  Moreover
coming back to inequality (\ref{cont3}), we deduce estimate
(\ref{uniqbis}), which ends the proof of the Lemma.
  \hfill $\Box$

\medskip

\subsection{Existence of solutions to problem  $({\mathcal P}_{per,n}^\varepsilon)$}
Let $\varepsilon >0$ and $n\geq 1$ be fixed. We shall prove  the
following result:
\begin{prop}\label{thm2}
Under  hypotheses (\ref{hyp}), (\ref{hyp zeta}) and (\ref{hyp chi}),
there exists a $T$-periodic solution to  problem $({\mathcal
P}_{per,n}^\varepsilon)$  denoted $(U_n^\varepsilon, \widetilde
\tau_n^\varepsilon, H_n^\varepsilon=\nabla \varphi_n^\varepsilon)$
such that
$$\begin{array}{c}U_n^\varepsilon\in {\mathcal C}([0,T]; {\mathcal U}),\ \
\widetilde \tau_n^\varepsilon\in \mathcal C([0,T]; L^2(D))\cap
L^2(0,T;
H^1_{\Gamma^-})\cap L^\infty(D_T), \\[2ex]  H_n^\varepsilon=\nabla \varphi_n\in {\mathcal
C}([0,T]; {\mathbb L}^2(D))\cap L^\infty(0,T; {\mathbb
L}^3(D)).\end{array}
$$
Moreover it holds
\begin{equation}\label{ineq1}
0\leq \widetilde \tau_n^\varepsilon+\zeta\leq \tau_\star \ \
\mbox{a.e. in } D_T.
\end{equation}
\end{prop}

For the proof, we consider for $(U_{0n}, \tau_0)\in \mathcal
V^n\times \mathcal W_0$,  the solution $(U_n, \widetilde\tau_n,$ $
H_n=\nabla \varphi_n)$ of problem $({\mathcal P}_{n}^\varepsilon
(U_{0n}, \tau_0))$ provided by Proposition \ref{prop2}  and we
define the map ${\mathcal K}: (U_{0n}, \tau_0)\mapsto (U_n(T),
\widetilde\tau_n (T))$.  We will apply the Schauder's fixed point
theorem to this map, following several steps. We start by looking
for a stable set under ${\mathcal K}$. In the sequel, we use the
notations and results of the previous subsection.
\bigskip

\noindent
$\bullet$ {\bf The subset  $\mathcal C_{R_\varepsilon}$ of ${\mathbb R}^n\times L^2(D)$.}\\
For any $R>0$, we consider the  closed bounded and convex subset of
${\mathbb R}^n\times L^2(D)$
\begin{equation}\label{subset}
\mathcal C_{R}=\Big\{(\alpha, \vartheta)\in {\mathbb R}^n\times
\mathcal W_0; \ \rho_0 \big(|\alpha|^2+ c_p\|\vartheta\|^2\big)\leq
R^2  \Big\},
 \end{equation}
 and  the map
 \begin{equation}\label{map}
{\mathcal K} : {\mathbb R}^n\times\mathcal W_0\to {\mathbb
R}^n\times\mathcal W_0,\ \ \  {\mathcal K}(\alpha_{0}^n, \tau _0)=
(\alpha^n(T), \widetilde\tau _n(T)),
\end{equation}
where $(U_n=\sum_{j=1}^n \alpha_j^n(t) A_j, \ \widetilde\tau_n,
H_n=\nabla \varphi_n)$ is the  solution of problem $({\mathcal
P}_n^\varepsilon (U_{0n},\tau_{0}))$ with  $U_{n0}= \sum_{j=0}^n
\alpha_{0j}^n A_j$. Recall that $\|U_{0n}\|= |\alpha_{0}^n|$ and
$\|U_{n}\|= |\alpha^n|$.

\medskip

 We have the result:
\begin{lem}\label{lem5}
There exists $R_\varepsilon>0$ independent of $n$ such that if
$(\alpha_{0}^n,  \tau _{0})$ belongs to $\mathcal C_{R_\varepsilon}
$ then  $(\alpha^n(T), \widetilde \tau _n(T))$ belongs to $\mathcal
C_{R_\varepsilon}$, that is ${\mathcal K}({\mathcal
C}_{R_\varepsilon})\subset {\mathcal C}_{R_\varepsilon}$ and
the following estimates are satisfied 
\begin{equation}\label{bound6}
\displaystyle \displaystyle \rho_0\big(\|U_n(t)\|^2+c_p\|\widetilde
\tau _n (t)\|^2\big)\leq R_\varepsilon^2,\ \ \ t\in [0,T],
\end{equation}
\begin{equation}\label{bound10}
\int_0^T e^{ \gamma t}(  \|\nabla U_n(t)\|^2+ \|\nabla \widetilde
\tau _n(t)\|^2) \, dt \leq \frac{ e^{\gamma T}}{\gamma}\
R_\varepsilon^2,
\end{equation}
where $\gamma >0$ is independent of $n$ and $\varepsilon$, (see
(\ref{gamma}) below). Moreover there exists a constant
$C({n,\varepsilon}, R_\varepsilon)>0$  depending  on
$R_\varepsilon$ such that
\begin{equation}\label{bound101}
 \| U_n(t)\|_\infty \leq C({n,\varepsilon}, R_\varepsilon), \ \ t\in [0,T]. \end{equation}
\end{lem}
\proof
 Here, we need
 to refine the estimates obtained previously.  We start with the following equality  satisfied by $\widetilde \tau_n$
$$
\displaystyle\frac{\rho_0 c_p}{2}\frac{d}{dt}\|\widetilde \tau
_n\|^2+ \eta \|\nabla \widetilde \tau _n\|^2= -\rho_0c_p   \int_D
(U_n\cdot \nabla \zeta) \, \widetilde \tau _n\, dx+\int_D Z(\zeta)
\, \widetilde \tau _n\, dx, $$ and we use the $L^\infty$ bound given
in (\ref{b1}) to
   get $$
\displaystyle\frac{\rho_0 c_p}{2}\frac{d}{dt}\|\widetilde \tau_n
\|^2+ \eta \|\nabla \widetilde \tau_n \|^2\leq C(1+\|U_n\|), $$ so
by using Poincar\'e inequality, we obtain
\begin{equation}\label{ener1} \displaystyle\frac{\rho_0c_p}{2}\frac{d}{dt}\|\widetilde \tau_n \|^2+ \eta \|\nabla \widetilde \tau_n \|^2\leq
C+\frac{\mu}{4}\|\nabla U_n\|^2. \end{equation} Adding  estimates
(\ref{ener1}) and  (\ref{b2}) leads to
$$\displaystyle  {\rho_0}\frac{d}{dt} \|U_n\|^2+{\rho_0
c_p}\frac{d}{dt}\|\widetilde \tau_n \|^2+ \frac{\mu}{2} \|\nabla
U_n\|^2+2 \eta \|\nabla \widetilde \tau_n \|^2\leq C_\varepsilon ,
$$  and applying again Poincar\'e inequality  to write $\|U_n\|\leq
C_p\|\nabla U_n\|$ and $\|\widetilde \tau_n\|\leq C_p\|\nabla
\widetilde \tau_n\|$, we obtain
\begin{equation}
\begin{array}{rr}
\displaystyle\rho_0\frac{d}{dt} (\|U_n\|^2 +c_p\|\widetilde \tau
_n\|^2)+ \frac{\mu}{4\rho_0 C_p^2} ( \rho_0\| U_n\|^2) +
\frac{\eta}{\rho_0c_p
C_p^2} (\rho_0c_p\| \widetilde \tau _n\|^2)
+\frac{\mu}{4} \|\nabla U_n\|^2+ \eta
\|\nabla \widetilde \tau _n\|^2 \leq C_\varepsilon.
\end{array}\nonumber
\end{equation}
Hence setting \begin{equation}\label{gamma}\gamma=
\min\Big(\frac{\mu}{4\rho_0 C_p^2},\ \frac{\eta}{\rho_0c_p C_p^2},\
\frac{\mu}{4},\ \eta \Big),
\end{equation}
we deduce that
\begin{equation}
\begin{array}{r}
\displaystyle\frac{d}{dt} (\rho_0 \|U_n\|^2 +\rho_0 c_p\|\widetilde
\tau
_n\|^2)+\gamma (\rho_0\| U_n\|^2+\rho_0c_p \| \widetilde \tau _n\|^2) + 
\gamma(
\|\nabla U_n\|^2+ \|\nabla \widetilde \tau _n\|^2) \leq
{C_\varepsilon}.
\end{array}\nonumber
\end{equation}
The above inequality can be rewritten in the form
\begin{equation}
 \displaystyle  \frac{d}{dt}\Big(e^{ \gamma t}\rho_0 (\|{ U}_n(t)\|^2+  c_p\|\widetilde \tau _n (t)\|^2 )\Big)+ \gamma e^{ \gamma t}\big( \|\nabla U_n\|^2+
  \|\nabla \widetilde \tau _n\|^2\big) \leq C_\varepsilon\, e^{\gamma  t},\nonumber
\end{equation}
and an integration  between $0$ and $t$ yields to
\begin{equation}\label{ener2}
\begin{array}{rr}
\displaystyle e^{\gamma  t}\rho_0 (\|{ U}_n(t)\|^2+  c_p\|\widetilde
\tau _n (t)\|^2 \big)+\gamma \int_0^t e^{\gamma s}\big( \|\nabla
U_n(s)\|^2
+\\[2ex] \displaystyle
\|\nabla \widetilde \tau _n(s)\|^2\big) \, ds \leq 
\rho_0 (\|{ U}_{0n}\|^2+  c_p\|  \tau _{0}\|^2) +
\frac{C_\varepsilon}{\gamma} (e^{\gamma t}-1).\end{array}
\end{equation}
  In particular for $t=T$, we get
\begin{equation}\label{ener3}
\begin{array}{r}
\displaystyle \rho_0 (|{ \alpha}^n(T)|^2+ c_p\|\widetilde \tau _n
(T)\|^2) \leq e^{-{\gamma} T}\rho_0 \big(|{ \alpha}_{0}^n|^2+ c_p\|
\tau _{0}\|^2
\big) + 
\frac{C_\varepsilon}{\gamma} (1-e^{- {\gamma}T}).
\end{array}
\end{equation}
Let $R_\varepsilon>0$ such that $ R_\varepsilon ^2\geq \displaystyle
\frac{C_\varepsilon}{\gamma}$ and assume that  $(\alpha_{0}^n, \tau
_0) \in {\mathcal C}_{R_\varepsilon}$. Then (\ref{ener3}) implies
\begin{equation}
\begin{array}{ll}
\displaystyle \rho_0 (|\alpha^n(T)|^2+  c_p\|\widetilde \tau _n
(T)\|^2) \leq e^{-\gamma T}R_\varepsilon^2 +
\frac{C_\varepsilon}{\gamma} (1-e^{-\gamma T}) \leq R_\varepsilon^2,
\nonumber
\end{array}
\end{equation}
thus the stability of ${\mathcal C}_{R_\varepsilon}$ by the map
${\mathcal K}$ is proved. Moreover, coming back to  (\ref{ener2}) we
deduce that for all $t\in [0,T]$
\begin{equation}
\displaystyle \displaystyle \rho_0 (\|{ U}_n(t)\|^2+
c_p\|\widetilde \tau _n (t)\|^2)\leq R_\varepsilon^2\, e^{- \gamma
t} + \frac{C_\varepsilon}{\gamma} (1-e^{- \gamma t}) \leq
R_\varepsilon^2
\end{equation}
and
$$\gamma \int_0^T e^{\gamma  t}(  \|\nabla U_n(t)\|^2+
\|\nabla \widetilde \tau _n(t)\|^2) \, dt\leq   R_\varepsilon^2
+\frac{ C_\varepsilon}{\gamma} (e^{\gamma T}-1)\leq  R_\varepsilon^2
e^{\gamma T},$$ which leads to (\ref{bound10}). Finally, since
$\rho_0\|U_{0n}\|^2\leq R_\varepsilon ^2 $ then
 we get (\ref{bound101}) by using estimate (\ref{boundlinfty}), see proof of Lemma \ref{lem4}.
 This ends the proof of the lemma. \hfill $\Box$

\medskip
\noindent
$\bullet$ {\bf Continuity and compactness of the map ${\mathcal K}$. }\\
We consider the set $\mathcal C_{R_\varepsilon} $ defined by
(\ref{subset}), where $R_\varepsilon >0$ is provided by Lemma
\ref{lem5}. We will prove that
\begin{lem}\label{lem6} The map ${\mathcal K}$ defined by
(\ref{map}) is continuous   and compact on  $\mathcal
C_{R_\varepsilon} $.
\end{lem}
  \proof
We will drop the index $n$ for simplicity, unless it is necessary to specify it.\\
 Let for  $k=1,2$,
$(\alpha_{0}^{ k},\tau_{0}^k)\in \mathcal C_{R_\varepsilon}$,
  $U_{0 }^k=\sum_{j=1}^n \alpha_{0j}^{ k}A_j $
 and let  $(U^k =\sum_{j=1}^n \alpha_{j}^{ k}A_j,$ $  \ \widetilde\tau^k, H^k=\nabla \varphi ^k)$ the
solution of problem $({\mathcal P}_n^\varepsilon (U_{0
}^k,\tau_{0}^k))$ provided by Proposition \ref{prop2}. We  apply
inequality (\ref{uniq}) for  $t=T$ to get using estimate
(\ref{bound101})
$$\begin{array}{r}
 \|U ^1(T)-U^2(T)\|^2+  \|\widetilde\tau ^1(T)-\widetilde\tau ^2(T)\|^2 
  \leq
   Ce^{T\, C({n,\varepsilon}, R_\varepsilon)}\, \big(\|U_{0 }^1-U_{0 }^2\|^2 +
   \|\tau_{0}^1-\tau_{0}^2\|^2 \big),
\end{array}
$$
 that is
$$|\alpha^{ 1}(T)-\alpha^{ 2}(T)|^2+  \|\widetilde\tau^1(T)-\widetilde\tau ^2(T)\|^2 
  \leq     Ce^{T\, C({n,\varepsilon}, R_\varepsilon)}\, \big(|\alpha_{0}^{ 1}-\alpha_{0}^{ 2}|^2+
   \|\tau_{0}^1-\tau_{0}^2\|^2 \big).$$
  This  means that $\mathcal K$ is Lipschitz continuous on $\mathcal
C_{R_\varepsilon}$.  Let us prove the compactness. \\
Let $(\alpha_{0}^k, \tau_0^ k)_k\subset \mathcal C_{R_\varepsilon}$
 and let for all $k\in \mathbb N$, $(U^k=
\sum_{j=1}^n\alpha^k_jA_j,\ \widetilde\tau^k,\ H^k=\nabla
\varphi^k)$ be the solution of problem (${\mathcal
P}^\varepsilon_{n}(U_{0}^k,\tau_0^k) $)  where
$U_{0}^k=\sum_{j=1}^n\alpha^k_{0j} A_j$ . Let $(\tau_0^{k_m})$ be a
subsequence of $(\tau_0^k)_k$ such that
$$\tau_0^{k_m} \rightharpoonup \tau_0 \ \ \ \mbox{ weakly in }
L^2(D).$$ In the sequel, we will write simply $m$ instead of $k_m$.
Let $\mathcal K(\alpha_{0}^m, \tau_0^m)= (\alpha^m(T),
\widetilde\tau^m(T)) $, we shall prove that $(\alpha^m(T),
\widetilde\tau^m(T))_m $ admits
a  subsequence converging in $\mathbb R^n\times L^2(D)$. \\
 First, since $(\alpha^m(T), \widetilde\tau^m(T))_m\subset  {\mathcal
C}_{R_\varepsilon}$ then by using Bolzano-Weierstrass Theorem we can
extract a subsequence of
 $(\widetilde\alpha^m(T))$ not relabeled which converges in ${\mathbb R}^n$.\\
  Next, according to estimates (\ref{bound6}) and (\ref{bound10}),  the
sequence  $ (\widetilde\tau^m)_m$ is bounded  in   $L^\infty(0,T;
L^2(D) )\, \cap\, L^2(0,T; H^1_{\Gamma^-})$ and
  $(U^m)_m$ is  bounded in $L^\infty(0,T; \mathcal U_0)\, \cap \,  L^2(0,T; \mathcal U)$.
   Hence, for  subsequences not relabeled,  we have
$$
\begin{array}{l}
U^m\rightharpoonup U\  {\mbox weakly }-\star\  \mbox{ in }\, L^\infty(0,T;\mathcal U_0) \ \mbox{ and weakly in }\ L^2(0,T; {\mathcal U}),\\[2ex]
\widetilde \tau^m\rightharpoonup \widetilde \tau \   {\mbox weakly
}-\star\ \mbox{ in }\, L^\infty(0,T; L^2(D))\  \mbox{ and weakly in
}\ L^2(0,T; H^1_{\Gamma^-}).
 \end{array}
$$
Now since
$$\rho_0c_p\partial_t \widetilde \tau^m =   \nabla \cdot (-\rho_0c_p \ \widetilde \tau^m  U^m  +\eta  \nabla\widetilde \tau^m )-\rho_0c_p
 U^m \cdot \nabla \zeta,
  $$
    and due to the
$L^\infty$-bound of $\widetilde \tau^m$ given in (\ref{b1}),  $
(\widetilde \tau^m U^m )_m $ is   bounded   in $L^\infty(0,T;\mathbb
L^2(D))$ and then $(\partial_t \widetilde \tau^m )_m$ is   bounded
in $L^2(0,T; H^{-1}_{\Gamma^-})$.  Hence for a subsequence
$$
\partial_t \widetilde \tau^m \rightharpoonup \partial_t \widetilde \tau  \ \  \ \mbox{ weakly in } {L^2(0,T;
H^{-1}_{\Gamma^-})},
$$ and  using Aubin-Lions compactness lemma we deduce that
$$  \widetilde \tau^m \rightarrow \widetilde \tau\,\,
\mbox{strongly in}\,\, L^2(0,T; L^2(D)).
$$
  Then it is easy to pass to the limit as $m\to
\infty$ in the weak formulation of the   equation of $\widetilde
\tau^m$ and we get that $\widetilde \tau $ satisfies for all
$\theta\in{\mathcal C}^1([0,T]\times \overline D)$ with $\theta=0$
on $\Gamma_T^-$ and $ \theta(T)=0$ in $D$
 $$ \begin{array}{cl}
\displaystyle -\rho_0c_p  \int_ {D_T} \widetilde \tau
\partial_t \theta\, dxdt -  \rho_0c_p  \int_{D_T} \widetilde \tau    U    \cdot \nabla\theta\, dxdt +\eta  \int_{D_T} \nabla \widetilde \tau   \cdot \nabla \theta\, dxdt  \\[3ex]
\displaystyle =\rho_0c_p\int_D  \tau_0 \theta(0)\, dx-\rho_0c_p\int_
{D_T} (U \cdot \nabla \zeta)\, \theta\, dxdt+\int_ {D_T} Z(\zeta)\,
\theta\, dxdt.
\end{array}
$$ Moreover  since $\partial_t {\widetilde \tau  }\in L ^2(0,T; H ^{-1}_{\Gamma^-}
 )$, we deduce that $\widetilde \tau  \in \mathcal C([0,T]; L^
2(D))$ and  $\widetilde \tau$ satisfies the temperature equation with the initial condition $\widetilde \tau(0)=  \tau_0$ .\\
We introduce the truncation function $\psi \in {\mathcal
C}^1([0,T])$ defined by $\psi(t)= 0$ for $t\in [0,\delta ]$ and
$\psi(t)= 1$ for $t\in [T-\delta, T]$ where $0<\delta< \frac{T}{2}$.
We see that  $w^m =\psi(t)(\widetilde \tau^m-\widetilde \tau)$
satisfies the problem
\begin{equation}\label{tempbis}
\begin{array}{l}
\rho_0 c_p \partial_t w^m -\eta \Delta w^m= G^m \ \  \mbox{ in } D_T,\\[2ex]
 w^m=0\ \  \mbox{on}\ \ \Gamma^-_T, \ \  \eta \nabla w^m\cdot \nu =0
\ \mbox{on}\ \Sigma_T\cup \Gamma^+_T ,\\[2ex]
w^m(0)= 0 \ \  \mbox{ in } D,
\end{array}
\end{equation}
where \  $ G^m= -\rho_0 c_p \psi(t) \big(U^m\cdot \nabla \widetilde
\tau^m -U\cdot \nabla \widetilde \tau+  (U^m-U)\cdot\nabla
\zeta\big) + \rho_0c_p \psi^\prime(t) (\widetilde \tau^m-\widetilde
\tau) $ \
 is bounded in $L^2(0,T;
 L^2(D))$ since   by estimate (\ref{bound101}),  $(U^m)_m$ is
bounded in $ {\mathbb L}^\infty(D_T)$.\\
  The following equality satisfied by $w^m$
$$
\begin{array}{c}
\displaystyle   \frac{\rho_0c_p}{2}\|w^m(T)\|^2  +\eta \int_0^T \|\nabla w^m\|^2\, dt=   \int_{D_T} G^m(t)\,  w^m(t)\, dxdt,
\end{array}
$$implies that
$$\|\widetilde \tau^m(T)-\widetilde \tau(T)\|^2\leq C \, \|G^m\|_{L^2(D_T)}\|w^m\|_{L^2(D_T)},$$
 and knowing that
$$
 w^m \rightarrow 0\,\,
\mbox{strongly in}\,\,  L^2(D_T),
$$  we deduce that
$$
\widetilde \tau^m(T)-\widetilde \tau(T)\to 0\ \ \mbox{ strongly in
}\ \ L^2(D).
$$
Hence ${\mathcal K}: \mathcal C_{R_\varepsilon} \to \mathcal
C_{R_\varepsilon}$ is  compact.\hfill $\Box$

\bigskip
\noindent
$\bullet$  {\bf End of  proof of  Proposition \ref{thm2}.}\\
As a consequence of Lemma \ref{lem6} and   Schauder's fixed point
Theorem, we deduce that ${\mathcal K }$ has a fixed point denoted
$(\alpha_{0}^{n\varepsilon}, \tau_0^\varepsilon)$ in $\mathcal
C_{R_\varepsilon}$. Therefore denoting $(U_n^{\varepsilon}=
\sum_{j=1}^n \alpha_j^{n\varepsilon}\, A_j, \,
\widetilde\tau_n^{\varepsilon},\,
H_n^{\varepsilon}=\nabla\varphi_n^{\varepsilon})$ the solution of
problem $(\mathcal P_ n^\varepsilon (U_{0n}^{\varepsilon},
\tau_0^{\varepsilon}))$ where $U_{0n}^{\varepsilon}=\sum_{j=1}^n
\alpha_{0j}^{n\varepsilon}\, A_j $, we infer that
$$(\alpha_{0}^{n\varepsilon}, \tau_0^{\varepsilon})={\mathcal K
}(\alpha_{0}^{n\varepsilon}, \tau_0^{\varepsilon})=
(\alpha^{n\varepsilon}(T), \widetilde\tau_n^{\varepsilon}(T)), $$
which means that the solution $(U_n^{\varepsilon} ,
\widetilde\tau_n^{\varepsilon}, H_n^{\varepsilon})$ is $T-$
periodic. Moreover estimate (\ref{ineq1}) is satisfied as a
consequence of (\ref{b1}).\hfill $\Box$

\section{Solutions to problem $({\mathcal P}_{per})$: Proof of Proposition \ref{thm3}}
We consider the time-periodic solution $(U_n^\varepsilon,\widetilde
\tau_n^\varepsilon, H_n^\varepsilon=\nabla \varphi _n^\varepsilon)$
of $({\mathcal P}_{per,n}^\varepsilon)$  provided by Proposition
\ref{thm2}.  We shall establish bounds for
$(U_n^\varepsilon,\widetilde \tau_n^\varepsilon, H_n^\varepsilon)$
which are independent of $n$   allowing to pass to the limit as
$n\to +\infty$.

\bigskip
\noindent
$\bullet$ {\bf Uniform bounds of $(U_n^\varepsilon, \widetilde \tau_n^\varepsilon, H_n^\varepsilon)$.}\\
In addition to the $L^\infty$ estimate (\ref{ineq1}) satisfied by
$\widetilde \tau_n^\varepsilon$, we have:
\begin{lem}\label{lem7}
For all $n\geq 1$ and $\varepsilon>0$, it holds
\begin{equation}\label{bound7}
\begin{array}{ll}
\displaystyle \int_0^{T}\|\nabla U_n^\varepsilon(s)\|^2\, ds
+\int_0^{T}\|\nabla \widetilde \tau_n^\varepsilon(s)\|^2\, ds \leq
C,
\end{array}
\end{equation}
\begin{equation}\label{bound8}
\begin{array}{ll}
\displaystyle \|U_n^\varepsilon\|_{L^\infty(0,T; \mathcal U_0 )}+
\|\widetilde \tau_n^\varepsilon\|_{ L^\infty(0,T; L^2(D))}+
 \|H_n^\varepsilon\|_{L^\infty(0,T; \mathbb L^3(D))}\leq C,
\end{array}
\end{equation}
\begin{equation}\label{dertau}
\|\partial_t \widetilde \tau_n^\varepsilon\|_{L^2(0,T;
H^{-1}_{\Gamma^-})}\leq C,
\end{equation}
\begin{equation}\label{otimes}
 \|U_n^\varepsilon\otimes
U_n^\varepsilon\|_{L^\frac{4}{3}(0,T;\mathbb L^{2}(D))}\leq C.
\end{equation}
 \end{lem}
\proof $(U_n^\varepsilon,\widetilde \tau_n^\varepsilon,
H_n^\varepsilon)$ satisfies  equations
(\ref{systemn1})-(\ref{systemn2})-(\ref{magn}) where we replace the
initial conditions by the periodic ones
$U_n^\varepsilon(0)=U_n^\varepsilon(T),\widetilde \tau_n^\varepsilon(0)=\widetilde \tau_n^\varepsilon(T)$ .\\
We start with  the following equalities
 \begin{equation}\label{b1bis}
\begin{array}{c}
\displaystyle\frac{\rho_0}{2}\frac{d}{dt} \|U^\varepsilon_n\|^2+  {\mu} \|\nabla U^\varepsilon_n\|^2=\int_D {\mathcal S}^n_\varepsilon \cdot U^\varepsilon_n\, dx,\\[3ex]
\displaystyle\frac{\rho_0 c_p}{2}\frac{d}{dt}\|\widetilde
\tau^\varepsilon_n\|^2+ \eta \|\nabla \widetilde
\tau^\varepsilon_n\|^2= -\rho_0c_p   \int_D (U_n^\varepsilon\cdot
\nabla \zeta) \, \widetilde \tau _n\, dx+ 
\int_D Z(\zeta) \, \widetilde \tau _n\,
dx,\end{array}
\end{equation}
where $ {\mathcal S}_\varepsilon^{n}={\mathcal
S}_\varepsilon(\widetilde \tau_n^\varepsilon, |H_n^\varepsilon|)$ so
inequality (\ref{Sbound2}) leads to
$$\begin{array}{rr}
\displaystyle|\int_D {\mathcal S}^n_\varepsilon\cdot
U_n^\varepsilon\, dx|= \big|\langle {\mathcal S}^n_\varepsilon,
U_n^\varepsilon\rangle_{\mathcal U\, ^\prime \times \mathcal U}
\big| \leq C\, \big(1 +
 \|\nabla\widetilde \tau_n^\varepsilon\|\big) \, \|\nabla U_n^\varepsilon\|,
 \end{array}$$
  $C>0$ being independent of $n$ and $\varepsilon$. Using  the $L^\infty $ bounds (\ref{ineq1}) of  $\widetilde \tau^\varepsilon_n$ and  Poincar\'e inequality
  we get
  $$|- \rho_0c_p \int_D  (U^\varepsilon_n\cdot \nabla \zeta) \,
\widetilde \tau^\varepsilon_n\, dx +\int_D Z(\zeta) \, \widetilde
\tau _n\, dx|\leq C \|U^\varepsilon_n\|+C\leq C_\beta +\beta
\|\nabla U^\varepsilon_n\|^2,
$$  where $\beta  >0$
will be chosen later independently of $\varepsilon$ and $n$, and $C_\beta >0$ depends on 
 $\beta $ but not of $\varepsilon$ and $n$. Therefore
\begin{equation}\label{equ1}
\begin{array}{ll}
\displaystyle \frac{\rho_0}{2}\frac{d}{dt}\|U_n^\varepsilon(t)\|^2+
\frac{\mu}{2} \|\nabla U_n^\varepsilon(t)\|^2\leq C\, \big(1 +
 \|\nabla\widetilde \tau_n^\varepsilon\|^2\big),\\[3ex]
\displaystyle \frac{\rho_0c_p}{2}\frac{d}{dt}\|\widetilde
\tau_n^\varepsilon(t)\|^2+ \eta\, \|\nabla \widetilde
\tau_n^\varepsilon(t)\|^2 \leq C_\beta +\beta \|\nabla
U^\varepsilon_n\|^2.
\end{array}
\end{equation}
 Integrating inequalities (\ref{equ1}) between $0$ and $T$ and using
the periodicity of the solution, we get
\begin{equation}\label{equ2}
\begin{array}{cl}
\displaystyle \frac{ \mu }{2}\int_0^T  \|\nabla U_n^\varepsilon(s)\|^2\, ds  \leq   CT +  C\int_0^T  \|\nabla \widetilde \tau_n^\varepsilon(s)\|^2\, ds,  \\[3ex]
\displaystyle   \eta \int_0^T \|\nabla \widetilde
\tau_n^\varepsilon(s)\|^2\, ds \leq C_\beta  T + \beta  \int_0^T
\|\nabla U_n^\varepsilon(s)\|^2\, ds.
\end{array}
\end{equation}
 Combining the two inequalities of \eqref{equ2}, we obtain
\begin{equation}\label{equ21}
\displaystyle(\frac{ \eta \mu }{2}-C\beta ) \int_0^T  \|\nabla
U_n^\varepsilon(s)\|^2\, ds \leq \displaystyle CT(\eta+ C_\beta ),
\end{equation}
and letting $\beta =\displaystyle \frac{\eta \mu}{4C}$ leads to
(\ref{bound7}).

\noindent Now let us prove   estimate  (\ref{bound8}). Let $0\leq
t\leq T$,  adding  the two inequalities in \eqref{equ1} and letting
$\beta =\displaystyle \frac{\mu}{2}$ yield for $s\in [t, t+T]$
\begin{equation}\label{equ3}
\begin{array}{ll}
\displaystyle \frac{\rho_0}{2}\frac{d}{ds}
\big(\|U_n^\varepsilon(s)\|^2+c_p \|\widetilde
\tau_n^\varepsilon(s)\|^2 \big) \leq C+C_\beta  +C\|\nabla
\widetilde \tau_n^\varepsilon(s)\|^2.
\end{array}
\end{equation}
 Let $t<t_1< t+T$, we  integrate \eqref{equ3} between $t_1 $ and $t+T$ to obtain
\begin{equation}\begin{array}{r}
\displaystyle \frac{\rho_0}{2}(\|U_n^\varepsilon(t+T)\|^2+ c_p
\|\widetilde \tau_n^\varepsilon(t+T)\|^2) \leq
\frac{\rho_0}{2}(\|U_n^\varepsilon(t_1)\|^2+c_p
\|\widetilde \tau_n^\varepsilon(t_1) \|^2)+ \\[2ex] \displaystyle \int_{t}^{t+T} (C+C_\beta
+C\|\nabla \widetilde \tau_n^\varepsilon(s)\|^2)\,ds, \nonumber
\end{array}\end{equation} and the periodicity implies
\begin{equation}\begin{array}{r}
\displaystyle \frac{\rho_0}{2}(\|U_n^\varepsilon(t)\|^2+ c_p
\|\widetilde \tau_n^\varepsilon(t)\|^2) \leq
\frac{\rho_0}{2}(\|U_n^\varepsilon(t_1)\|^2+c_p
\|\widetilde \tau_n^\varepsilon(t_1) \|^2)+\\[2ex] \displaystyle  (C+C_\beta )T+C\int_0^T\|\nabla
\widetilde \tau_n^\varepsilon(s)\|^2\,ds. \nonumber
\end{array}\end{equation} Therefore inequality \eqref{bound7} and
Poincar\'e inequality allow to get for $ t_1>t$
\begin{equation}
\displaystyle  {\rho_0} (\|U_n^\varepsilon(t)\|^2+c_p \|\widetilde
\tau_n^\varepsilon(t)\|^2 )\leq  {\rho_0C_p} (\|\nabla
U_n^\varepsilon(t_1)\|^2+ c_p\|\nabla \widetilde
\tau_n^\varepsilon(t_1) \|^2)+ C,  \nonumber
\end{equation}
and integrating with respect to $t_1$ between $t$ and $t+T$ leads to
\begin{equation}
\begin{array}{r}
\displaystyle   {\rho_0T} \big(\|U_n^\varepsilon(t)\|^2+c_p
\|\widetilde \tau_n^\varepsilon(t)\|^2 \big)   \leq 
{\rho_0}C_p \int_t^{t+T}
\big(\|\nabla U_n^\varepsilon(t_1)\|^2+ c_p\|\nabla \widetilde
\tau_n^\varepsilon(t_1) \|^2 \big) \,dt_1+ CT, \nonumber
\end{array}
\end{equation}
so using   again the bound \eqref{bound7}, we get
\begin{equation}
\begin{array}{cl}
\displaystyle  \|U_n^\varepsilon(t)\|^2+ \|\widetilde
\tau_n^\varepsilon(t)\|^2 \leq \displaystyle K, \nonumber
\end{array}
\end{equation}
where $K>0$ is independent of $\varepsilon$ and $n$. Finally using
the bound  (\ref{est magbis}) to estimate  $H_n^\varepsilon$, we
obtain inequality (\ref{bound8}).\\
Now,  since
$$\rho_0c_p\partial_t \widetilde \tau_n^\varepsilon =   \nabla \cdot (-\rho_0c_p \ \widetilde \tau_n^\varepsilon  U_n^\varepsilon  +\eta  \nabla\widetilde \tau_n^\varepsilon )-\rho_0c_p
 U_n^\varepsilon \cdot \nabla \zeta+Z(\zeta),
  $$
 and  $ \widetilde \tau_n^\varepsilon  U_n^\varepsilon  $ is uniformly bounded with
respect to $n$  and $\varepsilon $  in $L^\infty(0,T;\mathbb
L^2(D))$ we deduce that  $(\partial_t \widetilde \tau_n^\varepsilon
)_n$ is uniformly bounded with respect to $n$ and $\varepsilon$  in
$L^2(0,T; H^{-1}_{\Gamma^-})$. \\
Next, since $U_n^\varepsilon$ is uniformly bounded in $L^2(0,T;
\mathbb L^6(D))\cap L^\infty(0,T;\mathbb L^2(D)) $ with respect to
$n$ and $\varepsilon$, by an interpolation result, we deduce that
 $U_n^\varepsilon$ is uniformly bounded in $L^\frac{8}{3}(0,T;
\mathbb L^4(D))$ and so   $ U_n^\varepsilon\otimes U_n^\varepsilon$
is uniformly bounded with respect to $n$ and $\varepsilon$ in
$L^\frac{4}{3}(0,T; \mathbb L^2(D))$.\hfill $\Box$

\bigskip
\noindent
$\bullet$ {\bf Convergence results  as $n\to +\infty$. } \\
Let $\varepsilon>0$ be fixed. As a consequence of Lemma \ref{lem7},
we have:
 \begin{cor}\label{cor1} There exists subsequences still
labeled $(U_n^\varepsilon)_n$ and $(\widetilde
\tau_n^\varepsilon)_n$ such that when  $n\to +\infty$
$$
\begin{array}{l}
U_n^\varepsilon\rightharpoonup U^\varepsilon\,\, \mbox {weakly}-
\star\ \mbox{in} \,\, L^\infty(0,T; \mathcal U_0)
\ \  \mbox { and weakly  in} \,\, L^2(0,T; {\mathcal U}),\\[3ex]
\widetilde \tau_n^\varepsilon\rightharpoonup \widetilde
\tau^\varepsilon\,\, \mbox {weakly}- \star\  \mbox{in} \,\,
L^\infty(0,T; L^2(D))
\ \  \mbox { and weakly in} \,\, L^2(0,T; H^1_{\Gamma^-}),\\[3ex]
H_n^\varepsilon \rightharpoonup H^\varepsilon\,\, \mbox{weakly}-
\star\ \mbox{in}\,\, L^\infty(0,T; \mathbb L^3(D)),\\[2ex]
\partial_t \widetilde \tau_n^\varepsilon \rightharpoonup \partial_t \widetilde \tau^\varepsilon  \ \  \ \mbox{ weakly in } {L^2(0,T;
H^{-1}_{\Gamma^-})}.
\end{array}
$$
\end{cor}
In addition, we get the following strong convergences:
\begin{lem}\label{convfortes}
\begin{align}
& U_n^\varepsilon\rightarrow U^\varepsilon\,\,
\mbox{strongly in}\,\, L^2(0,T; \mathcal U_0), \label{convfort1}\\[2ex]
& \widetilde \tau_n^\varepsilon \rightarrow \widetilde
\tau^\varepsilon\,\, \mbox{strongly
in}\,\, L^2(D_T), \label{convfort2}\\[2ex] & H_n^\varepsilon=\mathcal
H(\widetilde \tau_n^\varepsilon)\rightarrow H^\varepsilon=\mathcal
H(\widetilde \tau^\varepsilon)\,\, \mbox{strongly in}\,\,  \mathbb
L^2(D_T), \label{convfort3}\\[2ex]
& \mathcal S_n^\varepsilon = \mathcal S_\varepsilon(\widetilde
\tau_n^\varepsilon, |H_n^\varepsilon|) \to\mathcal
S^\varepsilon=\mathcal S_\varepsilon(\widetilde \tau^\varepsilon,
|H^\varepsilon|)\ \ \mbox{strongly in }   \mathbb L^2(D_T).
\label{convfort4}
\end{align}

Moreover the limit $(U^\varepsilon, \widetilde \tau^\varepsilon,
H^\varepsilon)$ satisfies
\begin{align} & 0\leq \widetilde \tau^\varepsilon+\zeta\leq  \tau_\star
\ \ \mbox{  a.e. in } D_T, \label{bound1ter}\\[2ex]
 & \displaystyle \int_0^{T}\|\nabla U ^\varepsilon(s)\|^2\, ds
+\int_0^{T}\|\nabla \widetilde \tau ^\varepsilon(s)\|^2\, ds \leq C,
\label{bound77}\\[2ex]
& \displaystyle \|U ^\varepsilon\|_{L^\infty(0,T; \mathcal U_0 )}+
\|\widetilde \tau ^\varepsilon\|_{ L^\infty(0,T; L^2(D))}+
 \|H ^\varepsilon\|_{L^\infty(0,T; \mathbb L^3(D))}\leq C,
 \label{bound88}\\[2ex]
&  \|\partial_t \widetilde \tau^\varepsilon\|_{L^2(0,T;
H^{-1}_{\Gamma^-})}\leq C, \label{dertau2}\\[2ex]
& \|\mathcal S^\varepsilon\|_{L^2(0,T; \mathcal U^\prime)}\leq C.
\label{convfort4bis}\end{align}
\end{lem}
 \proof Following \cite{AH5}, see Lemma 3.5 and \cite{JLL} , we
 get the strong convergence (\ref{convfort1}). Next  using Aubin compactness
lemma we deduce (\ref{convfort2}) where $\widetilde
\tau^\varepsilon\in \mathcal C([0,T]; L^2(D))$ and satisfies $0\leq
\widetilde \tau^\varepsilon+\zeta\leq  \tau_\star$ {  a.e. in }
$D_T$.  Lemma \ref{lem3} leads to (\ref{convfort3}) where
 $H^\varepsilon\in  \mathcal C([0,T];
\mathbb L^2 )\, \cap \, L^\infty(0,T;  \mathbb L^3(D))$
 and satisfies $
\|H^\varepsilon\|_{L^\infty(0,T; \mathbb L^3(D))}\leq C$. Then we
use inequality (\ref{Scont}) to  get (\ref{convfort4}) and according
to (\ref{Sbound2}), $\mathcal S^\varepsilon$ satisfies the  estimate
(\ref{convfort4bis}). We end the proof of the lemma thanks to the
bounds (\ref{bound7}), (\ref{bound8}) and (\ref{dertau}). \hfill
$\Box$

\bigskip

\noindent
$\bullet$ {\bf Passage to the limit in problem $({\mathcal P}_{per, n}^\varepsilon)$ as $n\to +\infty$. } \\
Now, we pass to the limit as $n\to +\infty$ in the
weak formulation of problem $({\mathcal P}_{per, n}^\varepsilon)$.\\
First by Lemma \ref{lem3}, $H^\varepsilon=
\nabla\varphi^\varepsilon$ satisfies the magnetostatic equation
$$\int_D \Big(\nabla \varphi^\varepsilon+   M_S\, b(\widetilde \tau^\varepsilon(t)) \,
a(\nabla\varphi^\varepsilon)\Big) \cdot \nabla \phi \, dx= - \int_D
F(t)\, \phi\, dx, \;\, \forall \phi \in H^1_\sharp, \ \ t\in[0,T].$$
To deal with  equation of the velocity $U^\varepsilon$, we will use
the space $${\mathcal C}^1_{per}([0,T])= \{\psi \in {\mathcal C}^1
([0,T]); \psi(0)=\psi(T)\}.$$
 Let $\Psi\in {\mathcal U}$,
$\Psi_n\in \mathcal V^n$
 such that $\Psi_n\rightarrow \Psi$
strongly in ${\mathcal U}$ and let $\psi\in {\mathcal
C}^1_{per}([0,T])$,  we have for all $n\geq 1$
\begin{equation}\label{weakUn}
\begin{array}{ll}
\displaystyle -\rho_0 \int_{D_T} \psi^\prime(t)(
U_n^\varepsilon\cdot \Psi_n) \, dtdx-\rho_0 \int_{D_T}  \psi(t)(
U_n^\varepsilon\otimes U_n^\varepsilon\cdot \nabla \Psi_n)
 \, dtdx +\\[3ex]
\displaystyle  \mu\int_{D_T}\psi(t)( \nabla U_n^\varepsilon\cdot
\nabla \Psi_n) \, dtdx = \int_{D_T}\psi(t) ( {\mathcal
S}_n^\varepsilon\cdot \Psi_n) \, dtdx.
\end{array}
\end{equation}
The convergence (\ref{convfort1}) implies
$$U_n^\varepsilon\otimes  U_n^\varepsilon \rightarrow U^\varepsilon\otimes U^\varepsilon\,\, \mbox{strongly in}\,\, L^1(0,T;\mathbb  L^{1}(D)),$$
but    $ U_n^\varepsilon\otimes U_n^\varepsilon$ is uniformly
bounded with respect to $n$ and $\varepsilon$ in $L^\frac{4}{3}(0,T;
\mathbb L^2(D))$  by (\ref{otimes}) so we deduce that
\begin{equation}\label{convbis} U_n^\varepsilon\otimes
U_n^\varepsilon \rightharpoonup U^\varepsilon\otimes
U^\varepsilon\,\, \mbox{weakly in}\,\, L^\frac{4}{3}(0,T;\mathbb
L^{2}(D)), \end{equation} and
$$\|U^\varepsilon\otimes U^\varepsilon\|_{L^\frac{4}{3}(0,T;\mathbb
L^{2}(D))}\leq C.$$ Hence, passing to the limit as $n\to \infty$ by
using the previous convergence results, we obtain that
$U^\varepsilon$ verifies the equation
\begin{equation}\label{weakUeps}
\begin{array}{ll}
\displaystyle -\rho_0 \int_{D_T} \psi^\prime(t)( U ^\varepsilon\cdot
\Psi ) \, dtdx-\rho_0 \int_{D_T}  \psi(t)( U ^\varepsilon\otimes U
^\varepsilon\cdot \nabla \Psi )
 \, dtdx+ \\[3ex]
\displaystyle  \mu\int_{D_T}\psi(t)( \nabla U ^\varepsilon\cdot
\nabla \Psi ) \, dtdx = \int_{D_T}\psi(t) ( {\mathcal S}
^\varepsilon\cdot \Psi ) \, dtdx,\end{array}
\end{equation}
for all  $\Psi\in {\mathcal U}$  and  $\psi\in {\mathcal
C}^1_{per}([0,T])$. Integrating by parts with respect to $t$ the
first term, we deduce that $U^\varepsilon$ satisfies in ${\mathcal
D}^\prime(]0, T[)$
\begin{equation}\label{weakUeps2}
\begin{array}{r}
\displaystyle \rho_0\frac{d}{dt} \int_{D}{
U}^\varepsilon\cdot{\Psi}\, dx +\rho_0 \int_{D} ({
U^\varepsilon}\cdot\nabla){ U^\varepsilon} \cdot{\Psi}\, dx + \\[2ex] \displaystyle \mu
\int_{D} \nabla{ U^\varepsilon}\cdot \nabla {\Psi}\, dx  =  \int_D
{\mathcal S^\varepsilon}\cdot {\Psi}\, dx,
\end{array}
\end{equation} for all  $\Psi\in {\mathcal U}$, therefore
$\frac{d}{dt} \int_{D}{ U}^\varepsilon\cdot{\Psi}\, dx \in L^1(0,T)$
and so $\int_{D}{ U}^\varepsilon\cdot{\Psi}\, dx \in W^{1,
1}(0,T)\subset \mathcal C([0,T])$. Now we multiply equation
(\ref{weakUeps2}) by  $\psi\in {\mathcal C}^1_{per}([0,T])$   and
integrate between 0 and $T$, comparing the result with
(\ref{weakUeps}), we deduce that for all $\Psi\in {\mathcal U}$
\begin{equation}\label{perU}\int_{D}\big(U^\varepsilon(T)\psi(T)-U^\varepsilon(0)\psi(0)\big)\cdot
\Psi\, dx=0, \end{equation} so
$$U^\varepsilon(T)=U^\varepsilon(0).$$
Now we turn our attention to the temperature equation. Since $
\widetilde \tau_n^\varepsilon\, U_n^\varepsilon$ is uniformly
bounded with respect to $n$ and $\varepsilon$ in $L^2(0,T; \mathbb
L^2(D))$ and converges strongly towards $ \widetilde
\tau^\varepsilon\, U^\varepsilon$ in $L^1(0,T; \mathbb L^1(D))$ so
up to a subsequence
\begin{equation}\label{convter}
\widetilde \tau_n^\varepsilon\,U_n^\varepsilon\rightharpoonup
 \widetilde \tau^\varepsilon\, U^\varepsilon \,\, \mbox{weakly in}\,\, L^2(0,T; \mathbb
L^2(D)).
\end{equation}
Then it is easy to pass to the limit as $n\to \infty$ in the weak
formulation of the temperature equation and we get that $\widetilde
\tau^\varepsilon$ satisfies for all $\theta\in{\mathcal
C}^1([0,T]\times \overline D)$ with $\theta=0$ on $\Gamma_T^-$ and
$\theta(0)=\theta(T)$ in $D$
 \begin{equation}\label{weaktempeps}
\begin{array}{r}
\displaystyle -\rho_0c_p  \int_ {D_T} \big(\widetilde
\tau^\varepsilon
\partial_t \theta+ \widetilde \tau^\varepsilon  U^\varepsilon  \cdot \nabla\theta\big)\, dxdt + \eta  \int_{D_T} \nabla \widetilde \tau^\varepsilon \cdot \nabla \theta\, dxdt  \\[3ex]
\displaystyle = -\rho_0c_p\int_ {D_T}   (U^\varepsilon \cdot \nabla
\zeta)\, \theta\, dxdt+\int_ {D_T}   Z(\zeta)\, \theta\, dxdt.
\end{array}
\end{equation}
We deduce that $\widetilde \tau^\varepsilon$ satisfies in the sense
of distributions the equation
\begin{equation} \label{eqtauepsilon}\displaystyle \rho_0 c_p(\partial_t
\widetilde \tau^\varepsilon+ U^\varepsilon\cdot \nabla {\widetilde
\tau^\varepsilon}) -\eta \Delta \widetilde \tau^\varepsilon=
 -\rho_0 c_p \,  U^\varepsilon \cdot\nabla \zeta+Z(\zeta),
 \end{equation}
then   multiplying equation (\ref{eqtauepsilon}) by  function
$\theta\in{\mathcal C}^1([0,T]\times \overline D)$ with $\theta=0$
on $\Gamma_T^-$ and $\theta(0)=\theta(T)$ in $D$ we deduce the
relation
\begin{equation}
\int_D(\widetilde \tau^\varepsilon(T)\theta(T)- \widetilde
\tau^\varepsilon(0)\theta(0))\, dx=0,
\end{equation}
which implies the time-periodicity condition for $\widetilde
\tau^\varepsilon$ namely
\begin{equation}
\widetilde \tau^\varepsilon(T)= \widetilde \tau^\varepsilon(0),
\end{equation}
which in turn implies the time-periodicity condition for
$H^\varepsilon$.  We conclude that $(U^\varepsilon,\widetilde
\tau^\varepsilon, H^\varepsilon)$ is a weak solution of problem
($\mathcal P^\varepsilon_{per}$) satisfying   estimates
(\ref{bound9}).  To obtain estimate (\ref{bound9bis}) of $\partial_t
U^\varepsilon$, we use equation (\ref{weakUeps2}) and  the bounds of
$\mathcal S^\varepsilon$ in $ L^2(0,T; \mathcal U\, ^\prime)$ and
$U^\varepsilon\otimes U^\varepsilon$ in $ L^\frac{4}{3}(0,T;\mathbb
L^{2}(D))$. This ends proof of Proposition \ref{thm3}. \hfill $\Box$
\section{End of proof of Theorem \ref{thm1}}\label{endproof}
 Let $(U^\varepsilon, \widetilde \tau^\varepsilon, H^\varepsilon=\nabla
 \varphi^\varepsilon)$ be the solution of problem ($\mathcal
 P^\varepsilon_{per}$) provided by Proposition \ref{thm3}. The estimates
 obtained there allow to infer that up to subsequences not
 relabeled, we have  as $\varepsilon\to 0$
$$
\begin{array}{c}
(U^\varepsilon,\widetilde \tau^\varepsilon, H^\varepsilon) \rightharpoonup (U, \widetilde \tau, H)\,\, \mbox{weakly}-\star \ \mbox{in}\,  L^\infty(0,T;\mathcal U_0\times  L^2(D)\times \mathbb L^3(D)),\\[3ex]
 \widetilde \tau^\varepsilon\rightarrow  \widetilde \tau \,\, \mbox{strongly in} \,  L^2(0,T;
 L^2(D)),\\[3ex]
 \mathcal S^\varepsilon=\mathcal S_\varepsilon(\widetilde \tau^\varepsilon,
|H^\varepsilon|) \rightharpoonup  \mathcal S \,\, \mbox{ weakly in }
L^2(0,T;\mathcal U\, ^\prime).
\end{array}
$$
Moreover
$$0\leq \widetilde \tau +\zeta\leq  \tau_\star \ \  \mbox{ a.e. in }
D_T,$$ and by using Lemma \ref{lem3}, we deduce that
\begin{equation}
H^\varepsilon= \nabla \varphi^\varepsilon=\mathcal H(\widetilde
\tau^\varepsilon)\rightarrow H= \nabla \varphi =\mathcal
H(\widetilde \tau )\,\, \mbox{strongly in}\,\, L^2(0,T; \mathbb
L^2(D)),
\end{equation}
so  $\varphi$ satisfies (\ref{variat}). \\
 We will complete the passage to the limit as $\varepsilon \to 0$,
 once we prove the following results:
 \begin{lem}\label{derU}
 As  $\varepsilon \to 0$, we have
 \begin{align} & \label{convU}
U^\varepsilon  \rightarrow U\,\, \mbox{ strongly in }
L^2(0,T;\mathcal U_0), \\[2ex]
& \label{convUU} U^\varepsilon\otimes U^\varepsilon  \rightarrow
U\otimes U\,\, \mbox{ strongly in }
L^\frac{4}{3}(0,T;\mathbb L^2(D)),\\[2ex]
 & \label{convUtau} \widetilde \tau^\varepsilon U^\varepsilon\rightharpoonup \widetilde \tau
U \,\, \mbox{ weakly in }
  \mathbb L^2(D_T),\\[2ex]
& \label{convkappa} \kappa(\sigma ^\varepsilon \star|H^\varepsilon
|)\rightarrow \kappa(|H|) \,\, \mbox{strongly in} \
 L^2(D_T),\\[2ex]
 & \label{convS}
\mathcal S^\varepsilon=\mathcal S_\varepsilon(\widetilde
\tau^\varepsilon, |H^\varepsilon|) \rightharpoonup  \mathcal
S=\mathcal S (\widetilde \tau , |H| ) \,\, \mbox{ weakly in }
L^2(0,T;\mathcal U\, ^\prime).\end{align}
 \end{lem}
  \proof Since by Proposition \ref{thm3},  $U^\varepsilon$ is bounded in
$L^2(0,T; \mathcal U)$ and $\partial_t U^\varepsilon$ is bounded in
$L^\frac{4}{3}(0,T; \mathcal U\, ^\prime)$, then by Aubin
compactness lemma, we get the strong convergence (\ref{convU}). The
proofs of (\ref{convUU}) and (\ref{convUtau}) are  similar to those
of (\ref{convbis}) and (\ref{convter}). To  prove  convergence
(\ref{convkappa}), we write
$$
\begin{array}{r}\big|\kappa(\sigma ^\varepsilon \star|H^\varepsilon |)-
\kappa(|H|)\big|\leq \chi_0\big|\sigma ^\varepsilon
\star|H^\varepsilon |- |H|\, \big| \leq \\[2ex] \displaystyle
\chi_0\, \big|\sigma
^\varepsilon \star(|H^\varepsilon |- |H|)\big| +\chi_0\, \big|\sigma
^\varepsilon \star |H|-|H|\, \big| ,
\end{array}
$$ which leads to
$$\big\|\kappa(\sigma ^\varepsilon \star|H^\varepsilon |)-  \kappa(|H|)\big\|_{L^2(D_T)}\leq \chi_0\Big(\| H^\varepsilon - H\|_{L^2(D_T)}
+\big\|\sigma ^\varepsilon \star |H|-|H|\big\|_{L^2(D_T)}\Big).  $$
We know that the limit of the first term is 0 and to deal with the
second one, we see that for all $t\in [0,T]$
$$\displaystyle \lim_{\varepsilon \to 0}\|\sigma^ \varepsilon \star |H(t)|-|H(t)|\, \big\| =0 \ \
\mbox{ and }\ \  \|\sigma^ \varepsilon \star |H(t)|-|H(t)|\, \big\|
\leq 2 \|H(t)\| , \ \ \forall \varepsilon>0,$$ so by  Lebesgue's
dominated convergence Theorem $\|\sigma^ \varepsilon \star
|H(t)|-|H(t)|\, \big\| \to 0$ in $L^2(0,T)$, that is
$$\big\|\sigma^ \varepsilon \star |H(t)|-|H(t)|\, \big\|_{L^2(D_T)}\to 0.$$
 Now to prove (\ref{convS}), we
write, see the proof of estimate (\ref{Sbound2}) for details
$$
{\mathcal S}^\varepsilon =\mathcal S_\varepsilon(\widetilde
\tau^\varepsilon, |H^\varepsilon|) =\rho_0 (1+ \alpha\, b(\widetilde
\tau^\varepsilon))\, g+ \mu_0\, M_S \mathcal R^\varepsilon,
$$
where $$\mathcal R^\varepsilon   = \nabla \Big( b(\widetilde
\tau^\varepsilon)\kappa(\sigma ^\varepsilon \star|H^\varepsilon
|)\Big)-\kappa(\sigma ^\varepsilon \star|H^\varepsilon |) \nabla
b(\widetilde \tau^\varepsilon), $$ so for  $V\in L^2(0,T; \mathcal
U)$, we have
$$\begin{array}{r}  \displaystyle \int_{D_T} {\mathcal S^\varepsilon} \cdot V  \, dxdt= \rho_0\int_{D_T} (1+ \alpha\, b(\widetilde \tau^\varepsilon))\, g \cdot V  \, dxdt+
\\[2ex] \displaystyle \mu_0\, M_S\int_{D_T}
\kappa(\sigma ^\varepsilon \star|H^\varepsilon |) \nabla (\widetilde
\tau^\varepsilon + \zeta)\cdot V  \, dxdt.\end{array} $$ The results
above imply that as $\varepsilon \to 0$
\begin{equation}\label{convS1}\begin{array}{l} b(\widetilde \tau^\varepsilon)\rightarrow
b(\widetilde \tau) \,\, \mbox{strongly in} \ L^2(0,T;
 L^2(D)),\\[2ex]
 \nabla(\widetilde \tau^\varepsilon+\zeta)\rightharpoonup  \nabla(\widetilde \tau+\zeta) \,\, \mbox{weakly in} \  L^2(0,T;
 \mathbb L^2(D)).\\[2ex]
\end{array}\end{equation}
Moreover the inequalities  $\kappa(\sigma ^\varepsilon
\star|H^\varepsilon |)\leq \chi_0\  \sigma ^\varepsilon
\star|H^\varepsilon | $ \ and
$$\big\|\sigma ^\varepsilon \star|H^\varepsilon |\,
\|_{L^\infty(0,T;  L^3(D))}\leq \|H^\varepsilon\big\|_{L^\infty(0,T;
\mathbb L^3(D))},$$ imply that  $\kappa(\sigma ^\varepsilon
\star|H^\varepsilon |)$ is uniformly bounded in ${L^\infty(0,T;
L^3(D))}$. Therefore $\kappa(\sigma ^\varepsilon \star|H^\varepsilon
|) \nabla (\widetilde \tau^\varepsilon + \zeta)$ is uniformly
bounded in ${L^2(0,T; \mathbb L^\frac{6}{5}(D))}$ so up to a
subsequence $\kappa(\sigma ^\varepsilon \star|H^\varepsilon |)
\nabla (\widetilde \tau^\varepsilon + \zeta)$ converges weakly in $
L^2(0,T; \mathbb L^\frac{6}{5}(D))$. But (\ref{convkappa}) and
(\ref{convS1}) imply that $\kappa(\sigma ^\varepsilon
\star|H^\varepsilon |) \nabla (\widetilde \tau^\varepsilon + \zeta)$
converges towards $ \kappa(|H |) \nabla (\widetilde \tau  + \zeta) $
in the sense of distributions so we infer that
$$\kappa(\sigma ^\varepsilon \star|H^\varepsilon |) \nabla
(\widetilde \tau^\varepsilon + \zeta)\rightharpoonup \kappa(|H |)
\nabla (\widetilde \tau + \zeta)\ \ \mbox{weakly in }\  L^2(0,T;
\mathbb L^\frac{6}{5}(D)).
$$
Consequently for all $V\in L^2(0,T; \mathcal U)$
$$\int_{D_T} {\mathcal S^\varepsilon} \cdot V  \, dxdt\to  \rho_0\int_{D_T} (1+ \alpha\, b(\widetilde \tau ))\, g \cdot V  \, dxdt+\mu_0\, M_S\int_{D_T}
\kappa( |H  |) \nabla (\widetilde \tau  + \zeta)\cdot V  \, dxdt. $$
Now
$$\kappa( |H  |) \nabla (\widetilde \tau  + \zeta)= -\kappa( |H  |)\, \nabla (b(\widetilde \tau ))=-\nabla \Big(
b(\widetilde \tau )\kappa( |H  |)\Big)+b(\widetilde \tau)\chi(|H|)
\nabla |H|,
$$
so the previous limit reads as
$$\begin{array}{r}  \displaystyle \int_{D_T} {\mathcal S^\varepsilon} \cdot V  \, dxdt\to  \rho_0\int_{D_T} (1+ \alpha\, b(\widetilde \tau ))\, g \cdot V  \, dxdt+
\\[2ex] \displaystyle \mu_0\,
M_S\int_0^T\, \big \langle  b(\widetilde \tau)\chi(|H|) \nabla |H|,
V \big\rangle_{\mathcal U\, ^\prime \times \mathcal U}\,
dt,\end{array}$$ for all $V\in L^2(0,T; \mathcal U)$,  which is
nothing but (\ref{convS}).

\medskip

 Let us perform the limit as $\varepsilon \to 0$ in equation
 (\ref{weakUeps2}).   We multiply this equation  by $\psi
\in {\mathcal C} ^1_{per}([0,T])$   then we integrate by parts with
respect to $t$ and letting $\varepsilon\to 0$, we get using the
previous convergences
\begin{equation}\label{weakU} \begin{array}{ll} \displaystyle -\rho_0
\int_{D_T} \psi^\prime(t)( U  \cdot   \Psi ) \, dtdx-\rho_0
\int_{D_T} \psi(t)( U  \otimes U  \cdot \nabla \Psi )
 \, dtdx+ \\[3ex]
\displaystyle  \mu\int_{D_T}\psi(t)( \nabla U  \cdot \nabla \Psi )
\, dtdx = \int_0^T\, \psi(t)   \langle {\mathcal S},   \Psi
 \rangle
 \, dt,\end{array}
\end{equation}
 for all $\Psi \in \mathcal U$. Then proceeding as in the previous
section, we deduce that $U$ is $T$-periodic and satisfies
(\ref{eqU}).\\
Next it is easy to pass to the limit as $\varepsilon\to 0$ in
equation (\ref{weaktempeps}) then we get that $\widetilde \tau $ is
$T$-periodic and satisfies (\ref{energyequality}). Hence,  we
conclude that $(U, \widetilde \tau, H)$ is a weak solution of
problem
($\mathcal P_{per}$) in the sense of Definition \ref{definition}.\\
 It remains to  deal with the pressure.
  Since
$$\begin{array}{c}
\partial_t U\in
W^{-1,\infty}(0,T; \mathbb L^2(D)), \ \ (U\cdot\nabla)U=\nabla
\cdot(U\otimes U)\in L^\frac{4}{3}(0,T; \mathbb H^{-1}(D)),\\[2ex]
\Delta U \in L^2(0,T; \mathbb H^{-1}(D)), \ \ {\mathcal S}\in
L^\infty(0,T; \mathbb H^{-1}(D)),
\end{array}$$
  then $$\rho_0
(\partial_t { U}+ ({U}\cdot\nabla){ U})-\mu\Delta { U}  - {\mathcal
S}\in W^{-1,\infty}(0,T; \mathbb H^{-1}(D)).$$ By   De Rham's
Theorem, see \cite{Tar, TE}, we deduce that there exits  $p\in
W^{-1,\infty}(0,T; L^2(D))$ such that
 $$ \rho_0 (\partial_t { U}+ ({U}\cdot\nabla){ U})-\mu\Delta
{ U}  - {\mathcal S} =-\nabla p.$$   This ends  proof of Theorem
\ref{thm1}. \hfill $\Box$

\section{Other boundary conditions for the temperature}\label{other}
Our main result can be easily extended to various boundary conditions satisfied by the temperature. Two special cases are discussed below.\\

\noindent {\bf 1 -  Fluid heated from below and above}. We consider
in the model equations (\ref{systembis}) the temperature equation
with Dirichlet boundary conditions on $\Gamma_T^\pm$ and Neumann
boundary condition on $\Sigma_T$
$$
\begin{array}{ll}\displaystyle \rho_0 c_p(\partial_t {\tau}+ { U}\cdot\nabla{\tau}) -\eta \Delta \tau=
0 \ \ \mbox{ in } D_T,\\[3ex] \tau =\zeta^\pm\  \mbox{ on } \
\Gamma_T^\pm, \ \ \eta\nabla \tau\cdot \nu =0\ \  \mbox{ on } \
\Sigma_T,\\[2ex]
\tau(0)=\tau(T)\ \  \mbox{ in } \  D,
\end{array}
$$
 the velocity $U$ still satisfying the homogeneous boundary condition
$U=0$ on $(0,T)\times \partial D$. This problem will be labeled
($\mathcal P_{per}^\pm$). The corresponding  initial boundary
problem   was discussed in \cite{AH5}, considering $\zeta^+$ and
$\zeta^-$ are constants for simplicity. We assume that
$\widehat{\nabla}\zeta^\pm\cdot\widehat{\nu}=0$ on $(0,T)\times
\partial \Omega$   and $\zeta^\pm$  time-periodic with period $T$,
  and we  perform the change of variable  $\vartheta= \tau-\zeta$  where
$\displaystyle \zeta(t,x)= \frac{x_3}{d}\,
(\zeta^+(t,\widehat{x})-\zeta^-(t,\widehat{x}))+
\zeta^-(t,\widehat{x})$,  to get
 $$
\begin{array}{l}
\rho_0 c_p(\partial_t \vartheta+ U\cdot\nabla \vartheta)-\eta\
\Delta \vartheta =  - \rho_0 c_p U\cdot\nabla \zeta+Z(\zeta)
 \ \ \mbox{ in } D_T,\\[3ex] \vartheta =0\  \mbox{ on } \
\Gamma_T^\pm, \ \ \eta\nabla \vartheta\cdot \nu =0\ \  \mbox{ on } \
\Sigma_T,\\[2ex]
\vartheta(0)=\vartheta(T)\ \  \mbox{ in } \  D.
\end{array}
$$

Assuming that
\begin{equation}\label{hyp zeta bis}\begin{array}{ll}\zeta^\pm\in {\mathcal
C^1}([0,T];L^2(\Omega))\cap
\mathcal C([0,T]; H^2(\Omega) ),  \ \  \zeta^\pm \geq 0 \ \ \mbox{a.e. in } \Omega_T, \\[2ex] \displaystyle  \zeta^\pm(0)=\zeta^\pm(T)\ \ \mbox{ in }\Omega,
\ \ \widehat\nabla \zeta^\pm\cdot\widehat \nu=0 \ \mbox{ on }
(0,T)\times
\partial \Omega,
\\[2ex]
\tau_\star=\max(\|\zeta^-\|_{L^\infty(\Omega_T)},
\|\zeta^-\|_{L^\infty(\Omega_T)})>0,
\end{array}\end{equation}
we obtain an existence result    for this model which is similar to
that of Theorem \ref{thm1}.  More precisely, we have:
\begin{thm}\label{thm5}  Let $T>0$ be fixed. Under hypotheses (\ref{hyp}), (\ref{hyp chi}) and  (\ref{hyp zeta
bis}),  there exists a \, $T$-periodic solution $(U, \vartheta,
H=\nabla \varphi)$ of problem $({\mathcal P}_{per}^\pm)$ in the
sense  of  Definition \ref{definition} where we modify  the points
$(iii)$ and $(vii)$ as follows
\begin{enumerate}
\item[(iii)] the
temperature $\vartheta$ belongs to  ${\mathcal C}([0,T]; L^2(D))
\cap L^2(0,T; H^1_{\Gamma^\pm})$, $\vartheta(0)= \vartheta(T)$ and
satisfies $0\leq \vartheta+\zeta \leq \tau_\star\; \mbox{ a.e. in }
 D_T$, where $H^1_{\Gamma^\pm}=\left\{ \theta\in {H}^1(D); \,   \theta=0 \, \mbox{ on }
\Gamma^+\cup \Gamma^-\right\}$,
 \item[(vii)] the temperature $\vartheta$   satisfies
the integral identity
\begin{equation}
\begin{array}{lll}
\label{energyequality} && \displaystyle -\rho_0c_p\int_{D_T}
\vartheta \big(\partial_t \theta + { U} \cdot \nabla \theta\big) \,
dxdt + \eta \int_{D_T} \nabla \vartheta  \cdot \nabla \theta \, dxdt
\\[2ex]
&& \displaystyle = -\rho_0c_p  \int_{D_T} ({ U}\cdot \nabla \zeta)
\, \theta \, dxdt+ \int_{D_T}Z(\zeta) \, \theta \, dxdt,
\end{array}
\end{equation}
  for every ${\theta} \in
C^1([0, T] \times \overline{D})$ such that $\theta(0)=\theta(T)$ in
$D$ and $\theta=0$ on $\Gamma_T^-\cup \Gamma_T^+$.\end{enumerate}
Moreover there exists a pressure $p\in W^{-1,\infty}(0,T; L^2(D))$
such that the momentum equation is satisfied in the sense of
distributions.
\end{thm}
\proof The proof is similar to that of  Theorem \ref{thm1}. \hfill
$\Box$

\bigskip

\noindent {\bf 2 - Fourier boundary condition}. This time, we
consider the model equations (\ref{systembis}) with the temperature
satisfying  a Fourier boundary condition on $\Gamma_T^-$, that is
$\tau$ solves the following problem
$$
\begin{array}{ll}\displaystyle \rho_0 c_p(\partial_t {\tau}+ { U}\cdot\nabla{\tau}) -\eta \Delta \tau=
0 \ \ \mbox{ in } D_T,\\[3ex]  \  \tau+\eta\  \nabla \tau\cdot \nu = \zeta\  \mbox{ on } \
\Gamma_T^-, \ \ \eta\nabla \tau\cdot \nu =0\ \  \mbox{ on } \
\Sigma_T\cup \Gamma_T^+,\\[2ex]
\tau(0)=\tau(T)\ \  \mbox{ in } \  D,
\end{array}
$$
where $\zeta $ is time-periodic with period $T$. This new problem is labeled ($ \mathcal Q_{per}$). \\
Let us verify the energy estimate and  the maximum principle for the
corresponding initial boundary value problem for a given  initial
condition $\tau(0)= \tau_0$.
  The energy estimate derives from the equality
\begin{equation}\label{b111}
 \frac{\rho_0 c_p}{2}\frac{d}{dt}\|\tau\|^2+ \eta\|\nabla
\tau\|^2+ \int_\Omega|\tau(t,\widehat x, 0)|^2\, d\widehat{x}=
  \int_\Omega \zeta(t,\widehat x)\tau(t,\widehat x, 0) \,
d\widehat{x},
\end{equation}
which leads to
$$\begin{array}{r}
\displaystyle\rho_0 c_p  \|\tau(t)\|^2+ 2\eta\int_0^t\|\nabla
\tau(s)\|^2\, ds+
\int_0^t\|\tau(s, ., 0)\|_{L^2(\Omega)}^2\, ds \leq \\[2ex] \displaystyle  \rho_0 c_p
\|\tau_0\|^2+  \int_0^t\|\zeta(s,.)\|_{L^2(\Omega)}^2\, ds.
\end{array}$$ Next if $\zeta\geq 0$ in $\Omega_T$ and $0\leq \tau_0\leq
\tau_\star\equiv\|\zeta\|_{L^\infty(\Omega_T)}$ in $D$,
 then $\tau$
satisfies
$$
0\leq \tau \leq \tau_\star  \ \ \mbox{ a.e. in } D_T.
$$
 Hence assuming that
\begin{equation}\label{hyp zeta ter}\begin{array}{ll}\zeta \in {\mathcal C}([0,T];L^2(\Omega))  \cap L^\infty(\Omega_T),
 \ \  \zeta (0)=\zeta (T)\ \ \mbox{ in
}\Omega,  \\[2ex]  \zeta  \geq 0 \ \ \mbox{a.e. in } \Omega_T, \ \
\tau_\star= \|\zeta\|_{L^\infty(\Omega_T)} >0,
\end{array}\end{equation}
  we obtain the following existence result
  \begin{thm}\label{thm4}  Let $T>0$ be fixed. Assume    hypotheses (\ref{hyp}),  (\ref{hyp chi}) and (\ref{hyp zeta
ter})  to be hold true. Then there exists a \, $T$-periodic solution
$(U, \tau, H=\nabla \varphi)$ of problem $({\mathcal Q}_{per})$ in
the sense of  Definition \ref{definition} where instead of the
function $b$ defined in (\ref{b}), we take $b(\tau)=
\tau_\star-\tau$ and  we replace the points $(iii)$ and $(vii)$ by
the following ones
\begin{enumerate}
\item[$(iii)^\prime$] the temperature $\tau$ belongs to  ${\mathcal C}([0,T]; L^2(D)) \cap L^2(0,T; H^1(\Omega))$, $\tau(0)= \tau(T)$ and  satisfies $0\leq \tau \leq \tau_\star\; \mbox{ a.e. in }
 D_T$,
\item[$(vii)^\prime$] the temperature $\tau$   satisfies
the integral identity $$
\begin{array}{lll} && \displaystyle -\rho_0c_p\int_{D_T} \tau \big(\partial_t \theta
+ { U} \cdot \nabla \theta\big) \, dxdt + \eta \int_{D_T} \nabla
\tau  \cdot \nabla \theta \, dxdt
\\[2ex]
&& \displaystyle +  \int_{\Omega_T}  \tau(t,\widehat x, 0)\,
\theta(t,\widehat x, 0) \, d \widehat xdt =  \int_{\Omega_T}
\zeta(t,\widehat x)\, \theta(t,\widehat x, 0) \, d \widehat xdt,
\end{array}
$$
for any ${\theta} \in C^1([0, T] \times \overline{D})$ such that
$\theta(0)=\theta(T)$ in $D$.
\end{enumerate}

\noindent Moreover there exists a pressure $p\in W^{-1,\infty}(0,T;
L^2(D))$ such that the momentum equation is satisfied in the sense
of distributions.
\end{thm}
\proof  We adapt the proof of Theorem \ref{thm1} with the following
changes. The temperature equation has to be solved in the space
$H^1(D)$ and we will use the following inequality of Poincar\'e
type: there exists $C_p>0$ such that for all $w\in H^1(D)$ we have
$$
\|w\|^2\leq C_p^2\big( \|\nabla w\|^2 +\|w(., 0)\|_{L^2(\Omega)}^2\,
\big).
$$
 We replace the convex sets $\mathcal W_0$ and $\mathcal W$ by the
 following ones
\begin{equation}\label{tempspace2}
\begin{array}{c}\widehat{\mathcal W_0}=\big\{\vartheta \in L^\infty(D), \ \ \  0 \leq
\vartheta \leq \tau_\star\  \mbox{  a.e. in } D\big\},\\[2ex]
\widehat{\mathcal W}=\big\{\theta \in   \mathcal C([0,T];L^2(D)
)\cap L^\infty(D_T); \ \ \ 0\leq \theta \leq \tau_\star  \ \ \mbox{
a. e. in } D_T\big\},
\end{array}
\end{equation}
and we set $b(\tau)= \tau_\star-\tau$ so that if $\tau \in
\widehat{\mathcal W} $ then $b(\tau)\in \widehat{\mathcal W} $.
\\  With these new notations, all the results of the different steps
of proof of Theorem \ref{thm1} can be easily adapted to this case.
Nevertheless one must be more careful with the first part of proof
of Lemma \ref{lem5}, so for the
convenience of the reader, we will do it below.\\
Instead of $\widetilde \tau_n$ satisfying (\ref{systemn2}), we
consider  $\tau_n$  the solution of
\begin{equation}%
\begin{array}{ll}\displaystyle \rho_0 c_p(\partial_t {\tau_n}+ { U}\cdot\nabla{\tau_n}) -\eta \Delta \tau_n=
0 \ \ \mbox{ in } D_T,\\[3ex]  \  \tau_n+\eta\  \nabla \tau_n\cdot \nu = \zeta\  \mbox{ on } \
\Gamma_T^-, \ \ \eta\nabla \tau_n\cdot \nu =0\ \  \mbox{ on } \
\Sigma_T\cup \Gamma_T^+,\\[2ex]
\tau_n(0)=\tau_0\ \  \mbox{ in } \  D.
\end{array}
\end{equation} We use  equality (\ref{b111})  which leads to
\begin{equation}\label{ener111}
\displaystyle\frac{\rho_0 c_p}{2}\frac{d}{dt}\|\tau_n \|^2+ \eta
\|\nabla \tau_n \|^2+\frac{1 }{2} \|\tau_n(t,., 0)\|^2 \leq
\frac{1}{2 }\, \|\zeta(t)\|^2 \leq C.
\end{equation} Adding  estimates (\ref{ener111}) and (\ref{b2})
leads to
$$\displaystyle  {\rho_0}\frac{d}{dt} \|U_n\|^2+{\rho_0
c_p}\frac{d}{dt}\|\tau_n \|^2+ \frac{\mu}{2} \|\nabla U_n\|^2+2 \eta
\|\nabla \tau_n \|^2+  \|\tau_n(t,., 0)\|^2\leq C_\varepsilon ,
$$ and since $\|U_n\|\leq C_p\|\nabla U_n\|$ and $\|\tau_n\|^2\leq
C_p^2(\|\nabla \tau_n\|^2+\|\tau_n(t,., 0)\|_{L^2(\Omega)}^2)$, we
obtain
\begin{equation}
\begin{array}{rr}
\displaystyle\rho_0\frac{d}{dt} (\|U_n\|^2 +c_p\|\tau _n\|^2)+
\frac{\mu}{4\rho_0 C_p^2} ( \rho_0\| U_n\|^2) +
\frac{\eta_1}{\rho_0c_p
C_p^2} (\rho_0c_p\| \tau _n\|^2)\\[2ex]\displaystyle +\frac{\mu}{4} \|\nabla U_n\|^2+
\eta_2 \|\nabla \tau _n\|^2 +(1 -\eta_1) \|\tau_n(t,., 0)\|^2\leq
C_\varepsilon,
\end{array}\nonumber
\end{equation}
where  $ \eta_1, \eta_2>0$ are such that $ \eta_1+\eta_2=2\eta$.
Hence choosing $  0<\eta_1< \min( 1, 2\eta)$  and setting
\begin{equation}\label{gamma bis}\gamma= \min\Big(\frac{\mu}{4\rho_0
C_p^2},\ \frac{\eta_1}{\rho_0c_p C_p^2},\ \frac{\mu}{4},\ \eta_2
\Big),
\end{equation}
we deduce that
\begin{equation}
\begin{array}{r}
\displaystyle\frac{d}{dt} (\rho_0 \|U_n\|^2 +\rho_0 c_p\|\tau
_n\|^2)+\gamma (\rho_0\| U_n\|^2+\rho_0c_p \| \tau _n\|^2) +\\[2ex] \displaystyle
\gamma(\rho_0 \| U_n\|^2+\rho_0 c_p \|\tau _n\|^2) \leq
{C_\varepsilon}.
\end{array}\nonumber
\end{equation}
The  rest of the proof of Lemma \ref{lem5} remains valid.\\
We end by adapting the fundamental estimates of Lemma \ref{lem7} and
the
 $L^\infty$ estimate of the temperature
which read in the present case as follows
\begin{lem}\label{lem7bis} For all $n\geq 1$ and $\varepsilon>0$, it holds
\begin{equation}\label{bound14}
\begin{array}{ll}
0\leq \tau_n^\varepsilon \leq \tau_\star \ \ \mbox{a.e. in }
   D_T,
\end{array}
\end{equation}
\begin{equation}\label{bound15}
\begin{array}{ll}
\displaystyle \int_0^{T}\|\nabla U_n^\varepsilon(s)\|^2\, ds
+\int_0^{T}\big(\|\nabla
\tau_n^\varepsilon(s)\|^2+\|\tau^\varepsilon_n(s,.,0)\|^2_{L^2(\Omega)}\big)\,
ds \leq C,
\end{array}
\end{equation}
\begin{equation}\label{bound16}
\begin{array}{ll}
\displaystyle \|U_n^\varepsilon\|_{L^\infty(0,T; \mathcal U_0 )}+
\|\tau_n^\varepsilon\|_{ L^\infty(0,T; L^2(D))}+
 \|H_n^\varepsilon\|_{L^\infty(0,T; \mathbb L^3(D))}\leq C,
\end{array}
\end{equation}
\begin{equation}\label{bound17}
\|\partial_t \tau_n^\varepsilon\|_{L^2(0,T; (H^{1}(D))^\prime)}\leq
C,
\end{equation}
\begin{equation}\label{bound18}
 \|U_n^\varepsilon\otimes
U_n^\varepsilon\|_{L^\frac{4}{3}(0,T;\mathbb L^{2}(D))}\leq C.
\end{equation}
\end{lem}

\medskip

\noindent From here, to achieve the proof of Theorem \ref{thm4}, we
follow the same lines as in proving Theorem \ref{thm1}. \hfill
$\Box$

 \begin{rem} It would be interesting
to take into consideration other
 boundary conditions for the temperature and in particular non-linear
 conditions like Stephan-Boltzmann boundary condition. Likewise we can assume  that the velocity  satisfies the Navier boundary
 condition.
 \end{rem}

\section{Appendix}\label{appendix}
\noindent {\bf Proof of Proposition \ref{prop1}. }
  First, let us
give some useful properties satisfied by the function  $a( \xi)$
defined in (\ref{flux}).
 We make use of the
hypotheses (\ref{hyp chi}) for $\chi$  so  the function
$\xi\in{\mathbb R}^3\mapsto a(\xi)$ is continuous
 and we may write it as
\begin{equation}
a( \xi)= \nabla_\xi \Theta (\xi),\,\,  \Theta (\xi)= \kappa(|\xi|).
\end{equation}
Recall that  $\kappa$ is defined by (\ref{kappa}) and satisfies
(\ref{hyp kappa}) so $ \Theta\in \mathcal C^1(\mathbb R^3)$. Now
since $\kappa$  is an increasing function and the function
$\xi\mapsto |\xi|$ is convex, we deduce that
$\Theta$ is convex so we infer that function $ a $ is monotone. \\
The Hessian matrix of $\Theta$ denoted $h  = (h _{ij} )=
(\partial_{\xi_j}a_i )$  is given for all $i,j=1,2,3$ by
\begin{equation}
h _{ij}(\xi)=\left \{\begin{array}{ll}\displaystyle
\chi^\prime(|\xi|)\, \frac{\xi_i \xi_j}{|\xi|^2}+
\frac{\chi(|\xi|)}{|\xi|}\Big( \frac{\delta_{ij}|\xi|^2- \xi_i
\xi_j}{|\xi|^2}\Big) \ \ \mbox{if }\xi\neq 0,\\[3ex]
\chi^\prime(0) \, \delta_{ij} \ \ \mbox{if }\xi= 0,
\end{array}\right.
\end{equation}
where $\delta_{ij}$ is the Kronecker symbol.
 Setting $\xi =r\omega$ where
$r=|\xi|$ and $\omega\in S^1$ the unit sphere of $\mathbb R^3$, we
write
$$
h _{ij}(\xi)= \Big(\chi^\prime(r)-\frac{\chi(r)}{r}\Big)\
\omega_i\omega_j+\frac{\chi(r)}{r} \ \delta_{ij},
$$
and since  $ \displaystyle \lim _{r \to 0} \chi^\prime(r)= \lim _{r
\to 0}\frac{\chi(r)}{r}=\chi^\prime(0), $ \  we get\ $ \displaystyle
\lim _{\xi \to 0} h _{ij}(\xi)=  \chi^\prime(0) \, \delta_{ij}. $
Hence $\Theta \in \mathcal C^2(\mathbb R^3)$ and
 its  Hessian matrix $h  $ is bounded. Therefore $a \in \mathcal C^1(\mathbb R^3)$  and is Lipschitz continuous  on $\mathbb R^3$.\\
Next  we set for $t\in [0,T]$
\begin{equation}\label{atilde}
\widehat{a}(t,\xi)= \xi+  M_S\, b(\widetilde\tau(t, .)) \, a(\xi), \
\ \xi \in   {\mathbb R}^3,
\end{equation} then   $\widehat a(t,.)$ is Lipschitz continuous and  strongly monotone  on
${\mathbb R}^3$:
\begin{equation}
 \Big( \widehat a(t,\xi_1)-\widehat a(t,\xi_2)\Big) \cdot (\xi_1-\xi_2) \geq  |\xi_1-\xi_2|^2, \,\,
\forall\,\, \xi_1,\xi_2\in {\mathbb R}^3.
\end{equation}
\medskip

\noindent Now, we consider for $t\in [0,T]$ the weak formulation:
find $\varphi=\varphi(t)\in H^1_\sharp$ such that
\begin{equation}\label{magfaible}
 \int_D \Big(\nabla \varphi+   M_S\, b(\widetilde\tau(t)) \,
a(\nabla\varphi)\Big) \cdot \nabla \phi \, dx= - \int_D F(t)\,
\phi\, dx, \;\, \forall \phi \in H^1_\sharp.
\end{equation}
 We
set $$A(t, \varphi)\phi=\int_D \widehat a(t,\nabla \varphi) \cdot
\nabla \phi \, dx, \ \ \varphi, \phi\in H^1_\sharp, $$ where
$\widehat a(t, .)$ is given by (\ref{atilde}).  Then $A(t,.):
H^1_\sharp\to H^{-1}_\sharp$ is  continuous and  strongly monotone.
  So there exists a unique solution  $\varphi(t)\in H^1_\sharp $ of  (\ref{magfaible}) and we easily deduce that $\varphi(t)$ satisfies
  (\ref{magperlim}).   Now for
 $s, t\in [0,T]$,  using
the weak equations satisfied by $\varphi(s) $ and  $\varphi(t)$, we
get setting $K=\nabla \varphi(s)-\nabla\varphi(t) $
$$\begin{array}{c}
\|K\|^2\leq  \displaystyle \|K\|^2+ M_S\, \int_D   b(\widetilde
\tau(s)) \, \big( a(\nabla \varphi(s) )-a(\nabla \varphi(t))\big)
\cdot K\,
 dx= \\[2ex]\displaystyle M_S\,\int_D
(\widetilde \tau(t)-\widetilde \tau(s))\, a(\nabla \varphi(t)) \cdot
K\,
 dx - \int_D (F(s)-F(t))\, (\varphi(s)-\varphi(t))\, dx,
 \end{array}
$$
so the properties of  $a$  together with  Poincar\'e-Wirtinger
inequality allow to get
$$\|\nabla \varphi(s)-\nabla\varphi(t)\|\leq C (\|\widetilde \tau(t)-\widetilde \tau(s)\|+\|F(s)-F(t)\|).$$
Therefore $\varphi\in \mathcal C([0,T]; H^1_\sharp )$ and we easily
verify estimate (\ref{est mag}). Next writing equation
(\ref{magperlim}) as
$$
\Delta \varphi= F-\dvg G,
$$
where $G= M_S\, b(\widetilde \tau) \, a(\nabla\varphi)\in
L^\infty(D_T )$, we deduce using elliptic regularity results,  see
Lemma 4.27  in \cite{NS} and \cite{GM, GRIS} for example,  that
$\varphi\in L^\infty(0,T; W^{1,3}(D))$ and satisfies estimate
(\ref{est magbis}). \hfill $\Box$

\bigskip
\noindent {\bf Proof of Lemma \ref{lem3}. }  By using the
magnetostatic equations we have setting
$H_i=\mathcal H(\widetilde \tau_i)$
$$\begin{array}{r} \displaystyle \|H_1- H_2\|^2+ M_S\, \int_D   b(\widetilde \tau_1) ( a(H_1)-a(H_2)) \cdot
(H_1- H_2)\,
 dx= \\[2ex]\displaystyle  M_S\,\int_D
(\widetilde \tau_2-\widetilde \tau_1) a(H_2) \cdot (H_1- H_2)\,
 dx,\end{array}
$$
which leads to  inequality (\ref{Hcont}).  Therefore if $\widetilde
\tau^m \to \widetilde \tau \ \ \mbox{strongly
 in }
   L^2(D_T)$ and $\widetilde \tau \in \mathcal W$ then setting $H=\mathcal H(\widetilde \tau)$
   we see  that $H=\nabla \varphi$ where by
 Proposition \ref{prop1}, $\varphi  \in \mathcal C([0,T]; H^1_\sharp )\, \cap
\, L^\infty(0,T;  W^{1,3}(D))$,  solves equation (\ref{magfaible})
for all $t\in [0,T]$ and since by (\ref{Hcont})
$$\|H^m-H\|_{L^2(D_T)}\leq C \|\widetilde \tau^m-\widetilde \tau\|_{L^2(D_T)},$$
 then   $H^m  \to  H  \ \ \mbox{strongly in
   }
 \mathbb L^2(D_T)$. \hfill $\Box$

\bigskip
\noindent {\bf Proof of Proposition \ref{temp}.}     The  proof will
be done in two steps. First if we assume that
 $U\in   L^\infty(0,T; \mathcal U)$, these  results are classical, nevertheless for the reader's
convenience, we give a complete proof. We consider  for $t\in
(0,T)$, the bilinear form
 $$B(t; \widetilde \tau, \theta)= \int_D \Big(\rho_0 c_p ( U(t,.)\cdot\nabla{\widetilde \tau})\theta+\eta \nabla{\widetilde \tau}\cdot \nabla\theta \Big)\, dx, \ \ \
 \widetilde \tau, \theta \in H^1_{\Gamma^-},$$ and we set
$h(t)=  -\rho_0 c_p \ U(t,.)\cdot\nabla \zeta+Z(\zeta(t,.))$. Then
for every $\widetilde \tau, \theta \in H^1_{\Gamma^-}$ the function
$t\mapsto B(t; \widetilde \tau, \theta) $ is measurable and using
the embedding $H^1(D) \subset L^4(D)$ and Poincar\'e inequality we
get for a.e. $t\in (0,T)$
$$\begin{array}{rl}|B(t;   \widetilde\tau, \theta)| \leq \rho_0 c_p \|
U\|_{ L^\infty(0,T; L^4(D))} \ \|\nabla {  \widetilde\tau}\|\
\|\theta\|_{L^4(D)}+ \eta
\|\nabla{ \widetilde \tau}\| \| \nabla\theta \| \\[2ex]
\leq  C\| {
\widetilde\tau}\|_{H^1_{\Gamma^-}}\|\theta\|_{H^1_{\Gamma^-}}, \ \
\forall   \widetilde\tau, \theta \in H^1_{\Gamma^-},
\end{array}$$
where $C>0$ is a constant (depending upon $\| U\|_{ L^\infty(0,T;
L^4(D))}$) and
$$B(t; \widetilde \tau, \widetilde\tau)= \int_D \Big(\rho_0 c_p ( U(t,.)\cdot\nabla{\widetilde\tau}) \widetilde\tau +\eta |\nabla{\widetilde\tau}|^2 \Big)\, dx
=\eta \|\nabla{\widetilde\tau}\|^2, \ \ \forall   \widetilde\tau \in
H^1_{\Gamma^-} .$$ Since $h\in L^2(D_T)$,   applying Theorem of J.L.
Lions, see Theorem 10.9 in \cite{HB}, we get existence and
uniqueness of a weak solution $\widetilde \tau \in \mathcal C([0,T];
L^2(D))\cap L^2(0, T;H^1_{\Gamma^-})$,
$\displaystyle\frac{d\widetilde \tau}{dt}\in L^2(0,
T;H^{-1}_{\Gamma^-})$ satisfying
\begin{equation}\label{L1}\begin{array}{c} \displaystyle \rho_0 c_p\langle
\frac{d\widetilde \tau}{dt}, \theta\rangle+ B(t; \widetilde \tau,
\theta)=\langle  h(t), \theta\rangle,
 \ \ \mbox{ a.e. } t\in (0,T), \ \ \forall \theta \in H^1_{\Gamma^-}, \\[2ex]
\widetilde \tau(0)= \tau_0,  \end{array}\end{equation} and we
deduce, in a classical way,  that $\widetilde \tau$ satisfies
problem (\ref{tempinit}).
Next testing  equation (\ref{L1}) by $\widetilde \tau$   leads to
$$
\displaystyle\frac{\rho_0 c_p}{2}\frac{d}{dt}\|\widetilde \tau \|^2+
\eta \|\nabla \widetilde \tau \|^2=- \rho_0 c_p \ \int_D({
U}\cdot\nabla \zeta)\ \widetilde \tau \, dx+\int_D Z(\zeta)\
\widetilde \tau \, dx.
$$
Since $-\int_D({ U}\cdot\nabla \zeta)\ \widetilde \tau \,
dx=\int_D({ U}\cdot\nabla \widetilde \tau)\ \zeta \, dx$ then using
Cauchy-Schwarz, Poincar\'e and Young inequalities we get
$$
\displaystyle\frac{\rho_0 c_p}{2}\frac{d}{dt}\|\widetilde \tau \|^2+
\eta \|\nabla \widetilde \tau \|^2\leq
C(\|U\|^2\|\zeta\|_{L^\infty(\Omega_T)}^2+
\|Z(\zeta)\|_{L^\infty(0,T; L^2(\Omega))} ^2)+
\frac{\eta}{2}\|\nabla \widetilde \tau\|^2.
$$
Therefore
$$
 \displaystyle{\rho_0
c_p}\frac{d}{dt}\|\widetilde \tau \|^2+ \eta \|\nabla \widetilde
\tau \|^2\leq C(1+\|U\|^2).
$$
We get that $\widetilde \tau$ satisfies estimate (\ref{ener1bis}), by integration with respect to time.\\
 Next we
consider the function  $w= \widetilde \tau+\zeta$ which satisfies
the problem
$$
\begin{array}{l}
\rho_0 c_p (\partial_t w+ U\cdot\nabla w)-\eta \Delta w=0\,\, \mbox{ in } \, D_T,\\[2ex]
\displaystyle w =\zeta\ \mbox{ on } \ \Gamma_T^-,\ \ \eta\nabla
w\cdot \nu =0\ \  \mbox{ on } \ \Sigma_T\cup \Gamma_T^+, \\[2ex]
\displaystyle {w}(0)=w_0:={ \tau}_0+\zeta \ \  \mbox{ in } \ D.
\end{array}
$$
 We multiply the equation by the negative part $w^-$ of $w$, integrate by parts and
use the fact that $\zeta$ is nonnegative   to cancel the boundary
term, we get
$$
\frac{\rho_0 c_p}{2}\frac{d}{dt}\|w^-\|^2+ \eta\|\nabla w^-\|^2=0.
$$
Therefore since $w_0\geq 0$, we conclude that $w^-=0$ that is $w\geq
0$ a.e. in $D_T$. Next using the equation satisfied by $
\tau_\star-w$ we get as previously that   $(\tau_\star-w)^-=0$ so
$w\leq \tau_\star$ a.e. in $D_T$.\\
Now, we consider the general case where $U\in L^2(0,T; \mathcal U)$,
then by regularization, we can approximate $U$ by a sequence
$(U^m)_m\subset  L^\infty(0,T; \mathcal U)$ such that
$$U^m \to U\ \ \ \mbox{ strongly  in } L^2(0,T;  \mathcal U).$$
Indeed consider a regularizing sequence $\eta_m=\eta_m(t)$  and set
$U^m= \eta_m \ast \overline U$ where $\overline U$ is the extension
of $U$ by 0 outside $[0,T]$. From the first step, we deduce that for
each $m$, there exists a unique weak solution $\tau^m \in \mathcal
C([0,T]; L^2(D))\cap L^2(0, T;H^1_{\Gamma^-})\cap L^\infty(D_T)$ ,
$\displaystyle\frac{d\tau^m}{dt}\in L^2(0, T;H^{-1}_{\Gamma^-})$ of
problem
$$\begin{array}{ll}
\displaystyle \rho_0 c_p(\partial_t {\tau^m}+ { U^m}\cdot\nabla{\tau^m}) -\eta \Delta \tau^m= -\rho_0 c_p \ { U^m}\cdot\nabla \zeta+Z(\zeta) \,\, \mbox{ in } \, D_T,\\[2ex]
\displaystyle \tau^m =0\, \mbox{ on } \, \Gamma_T^-,\,\, \eta\nabla
\tau^m\cdot \nu =0\,\, \mbox{ on } \, \Sigma_T\cup \Gamma_T^+, \\[2ex]
\displaystyle {\tau^m}(0)={\tau}_0 \,\, \mbox{ in } \, D,
\end{array}$$
 satisfying
$$
0 \leq\tau^m+\zeta\leq \tau_\star \mbox{ a.e. in} \
D_T,$$
and
 $$
 \|\tau^m (t)\|^2+ \int_0^t \|\nabla
\tau^m(s) \|^2\, ds \leq C(1+\|\tau_0 \|^2+ \int_0^t \|U^m(s) \|^2\,
ds ), \ \ t\in [0,T].$$
 Hence $\tau ^m$ is uniformly
bounded in $L^\infty(D_T)\cap L^\infty(0,T; L^2(D)) \cap L^2(0,T;
H^{1}_{\Gamma^-})$. Moreover writing
$$\begin{array}{c} \rho_0 c_p\partial_t {\tau^m}
=\nabla\cdot (-\rho_0 c_p { U^m}{\tau^m} -\eta \nabla \tau^m -\rho_0
c_p \ { U^m}\cdot  \zeta)+Z(\zeta)\end{array}$$ where $U^m$ is
uniformly bounded in $L^2(0,T;  \mathcal U)$, we see that
 $\partial_t {\tau^m}$ is uniformly bounded in  $L^2(0,T; H^{-1}_{\Gamma^-})$.
Therefore, there exists a function $\widetilde \tau $ such that up
to a subsequence not relabeled, we have
$${\tau^m} \rightharpoonup \widetilde\tau \mbox { weakly-$\star$ in } L^\infty(0,T; L^2(D)) \mbox { weakly in } L^2(0,T; H^{1}_{\Gamma^-}),  $$
$$\partial_t {\tau^m}\rightharpoonup \partial_t \widetilde\tau \mbox { weakly in } L^2(0,T; H^{-1}_{\Gamma^-}),$$
$${\tau^m} \to \widetilde\tau \mbox { strongly in } L^2(0,T; L^4(D)),$$
$$0 \leq\widetilde \tau+\zeta\leq \tau_\star
\mbox{ a.e. in} \ D_T,$$ the strong convergence being obtained using
Aubin-Lions compactness lemma.  Hence $\widetilde \tau \in \mathcal
C([0,T]; L^2(D))\cap L^2(0, T;H^1_{\Gamma^-})\cap L^\infty(D_T)$ and
satisfies estimate (\ref{linfty}).\\ 
 In order
to pass to the limit in the equation satisfied by $\tau^m$, let us
consider test functions $\psi=\psi(t) \in \mathcal D([0,T[)$,
$\theta \in H^1_{\Gamma^-}$ and write using an integration by parts
$$\begin{array}{c}
\displaystyle -\rho_0 c_p\int_{D_T}    \tau^m\theta \psi^\prime(t)\,
dxdt  + \int_{D_T}  \Big(\rho_0 c_p ( U^m
\cdot\nabla{\tau^m})\theta+\eta \nabla{\tau^m}\cdot
\nabla\theta \Big)\psi(t)\, dxdt=\\[2ex]\displaystyle  \rho_0 c_p\int_{D}   \tau_0 \theta
\psi(0)\, dx+\int_{D_T}\Big(-\rho_0 c_p \ U^m\cdot\nabla
\zeta+Z(\zeta )\Big) \theta \,  \psi(t)\, dxdt.
\end{array}$$
  Since $$\int_{D_T}   ( U^m \cdot\nabla{
\tau^m})\theta \psi(t) \, dx dt = -\int_{D_T}   ( U^m \cdot\nabla
\theta){ \tau^m}\psi(t) \, dx
 dt\to  -\int_{D_T}   ( U \cdot\nabla \theta){  \tau}\psi(t) \, dx
 dt$$ as $m\to \infty$ then we get for all $\psi=\psi(t) \in \mathcal D([0,T[)$,
$\theta \in H^1_{\Gamma^-}$
$$\begin{array}{c}
\displaystyle -\rho_0 c_p\int_{D_T}    \widetilde\tau \theta
\psi^\prime(t)\, dxdt  + \int_{D_T}  \Big(\rho_0 c_p ( U
\cdot\nabla{\widetilde\tau })\theta+\eta \nabla{\widetilde\tau
}\cdot
\nabla\theta \Big)\psi(t)\, dxdt=\\[2ex]\displaystyle  \rho_0 c_p\int_{D}   \tau_0 \theta
\psi(0)\, dx+\int_{D_T}\Big(-\rho_0 c_p \ U \cdot\nabla
\zeta+Z(\zeta )\Big) \theta \,  \psi(t)\, dxdt.
\end{array}$$
This implies that $\widetilde\tau$ satisfies
\begin{equation}\label{L4}\begin{array}{c} \displaystyle \rho_0 c_p\langle
\frac{d\widetilde \tau}{dt}, \vartheta\rangle+ \int_D \Big(\rho_0
c_p ( U(t,.)\cdot\nabla{\tau})\vartheta+\eta \nabla{\tau}\cdot
\nabla\vartheta \Big)\, dx=\\[2ex] \displaystyle \int_D \big(-\rho_0 c_p \
U(t,.)\cdot\nabla \zeta+Z(\zeta(t,.))\big) \vartheta\, dx,
 \ \ \mbox{ a.e. } t\in (0,T), \ \ \forall \vartheta \in H^1_{\Gamma^-}, \\[2ex]
\widetilde \tau(0)= \tau_0,  \end{array}\end{equation}  and the
estimate (\ref{ener1bis}) is obtained as in the first case. It
remains to prove uniqueness. Let $\widetilde \tau ^1$ and
$\widetilde \tau ^2$ two solutions having the properties given in
the proposition. Then $\vartheta=\widetilde\tau^1-\widetilde\tau^2$
satisfies the problem
$$\begin{array}{ll}
\displaystyle \rho_0 c_p(\partial_t {\vartheta}+ { U}\cdot\nabla{\vartheta}) -\eta \Delta \vartheta= 0\ \  \mbox{ in } \  D_T,\\[2ex]
\displaystyle \vartheta =0\  \mbox{ on } \  \Gamma_T^-,\ \
\eta\nabla
\vartheta\cdot {\nu}=0\ \  \mbox{ on } \  \Sigma_T\cup\Gamma_T^+, \\[2ex]
\displaystyle {\vartheta}(0)=0\ \  \mbox{ in } \  D,
\end{array}
$$
so
$$
\displaystyle\frac{\rho_0 c_p}{2}\frac{d}{dt}\|\vartheta \|^2=- \eta
\|\nabla \vartheta \|^2,
$$ so that $\|\vartheta(t) \|=0$ for all $t\in [0,T]$. This ends the
proof of the proposition.  \hfill $\Box$

 \bigskip

 \noindent {\bf Proof of Lemma \ref{taucont}.}    $\theta=\widetilde\tau^1-\widetilde\tau^2$  satisfies the
problem
$$\begin{array}{ll}
\displaystyle \rho_0 c_p(\partial_t {\theta}+ { U^1}\cdot\nabla{\theta}) -\eta \Delta \theta= -\rho_0 c_p \ (U^1-U^2)\cdot \nabla({\widetilde\tau^2}+  \zeta )\ \  \mbox{ in } \  D_T,\\[2ex]
\displaystyle \theta =0\  \mbox{ on } \  \Gamma_T^-,\ \  \eta\nabla
\theta\cdot {\nu}=0\ \  \mbox{ on } \  \Sigma_T\cup\Gamma_T^+, \\[2ex]
\displaystyle {\theta}(0)=\tau_0^1-\tau_0^2 \ \  \mbox{ in } \  D,
\end{array}
$$
and using estimate (\ref{linfty}) for $\widetilde\tau^2$ , we get
\begin{equation}\label{contau1}\begin{array}{r}
\displaystyle\frac{\rho_0 c_p}{2}\frac{d}{dt}\|\theta \|^2+ \eta
\|\nabla \theta \|^2= \rho_0 c_p \ \int_D((U^1-U^2)\cdot\nabla
\theta)\ ({\widetilde\tau^2}+ \zeta ) \, dx\leq \\[2ex] \displaystyle C\|U^1-U^2\|^2
+\frac{\eta}{2}\|\nabla \theta\|^2,\end{array}
\end{equation}
which leads to estimate (\ref{contau}).\hfill $\Box$

\bigskip

\noindent {\bf Proof of Lemma \ref{Sprop}.} (\ref{Sbound1}) is
obtained using the bounds of $\chi$, $\, b$ and (\ref{est
  mag}) which imply
$$\|{\mathcal
S}_\varepsilon(t) \|\leq C_\varepsilon \|H(t)\|+ C\|g(t)\|\leq
C_\varepsilon \|F\|_{L^\infty(0,T; L^2(D))}+ C.$$ To prove the
second estimate, we  set $ \widetilde \tau=\widetilde
\tau^1-\widetilde \tau^2$ and $ H=H^1-H^2$ and write
$$
\displaystyle {\mathcal S}_\varepsilon^1-{\mathcal S}_\varepsilon^2
=\mu_0\, M_S   \,   \sum_{i=1}^3\ \mathcal R^i- \rho_0 \alpha  \,
\widetilde \tau g,
$$ where
 $$\begin{array}{ll}
\mathcal R^1=   b( \widetilde \tau^1)\ \chi(\sigma^\varepsilon
\star| H^1|)\ \nabla(\sigma^\varepsilon \star (|
H^1|-| H^2|)),\\[3ex]
\displaystyle  \mathcal R^2= -(\widetilde \tau_1-\widetilde \tau_2)
\ \chi(\sigma^\varepsilon \star| H^1|) \ \nabla(\sigma^\varepsilon
\star| H^2|),
\\[3ex]\displaystyle  \mathcal R^3=  b( \widetilde \tau^2) \ \Big(\chi(\sigma^\varepsilon \star| H^1|)-\chi(\sigma^\varepsilon \star| H^2|)\Big)\
\nabla(\sigma^\varepsilon \star| H^2|).
\end{array}$$
 Using properties of $\chi$ and  (\ref{est mag}), we see  that
$$
\begin{array}{ll}
\displaystyle \| \mathcal R^1\|\leq C_\varepsilon    \| H\|, \,\,\,
\displaystyle \| \mathcal R^2\|\leq C_\varepsilon
  \| \widetilde \tau\|,\,\, \,
\displaystyle \| \mathcal R^3\|\leq C_\varepsilon  \|  H\|,
\end{array}
$$
so the result follows by using (\ref{Hcont}). \\
 To prove the third estimate, we write
$$
{\mathcal S}_\varepsilon  =\rho_0 (1+ \alpha\, b(\widetilde \tau))\,
g+ \mu_0\, M_S \mathcal R,
$$
where
\begin{align*}
\mathcal R& = b(\widetilde \tau) \chi(\sigma ^\varepsilon \star|H
|)\nabla (\sigma^\varepsilon\star |H |) =
b(\widetilde \tau) \nabla \big(\kappa(\sigma ^\varepsilon \star|H |)\big)\\[2ex]
&=\nabla \big( b(\widetilde \tau)\kappa(\sigma ^\varepsilon \star|H
|)\big)-\kappa(\sigma ^\varepsilon \star|H |) \nabla ( b(\widetilde
\tau)),
\end{align*}
so for  $V\in \mathcal U$, we have
$$\int_D {\mathcal R} \cdot V  \, dx= \int_D \kappa(\sigma ^\varepsilon
\star|H |) \nabla (\widetilde \tau + \zeta)\cdot V  \, dx. $$ Since
$\kappa(s)\leq \chi_0 s$\  for $s\geq 0$, then  using
 the bound
(\ref{est magbis}) of $ H$ in   $\mathbb L^3(D)$ and the embedding
$\mathbb H^1 (D)\subset \mathbb L^6(D) $, we obtain
$$
\begin{array}{rr}
\displaystyle |\int_D {\mathcal S}_\varepsilon\cdot V  \, dx|\leq
C\,  \|V
\| + C\,  (\|\nabla\widetilde \tau  \| +\|\nabla \zeta  \|)\,  \|H  \|_3\, \|V  \|_6 \\[2ex]
\displaystyle \leq  C\, \big(1 +  \|\nabla\widetilde \tau   \| \big)
\|\nabla V \|,
\end{array}
$$
which ends the proof of the lemma. \hfill $\Box$




\end{document}